\documentclass{article}%
\usepackage{amsfonts}
\usepackage{amsmath}
\usepackage{amssymb}
\usepackage{amsxtra}
\usepackage{graphicx}
\usepackage{geometry}
\usepackage{color}
\usepackage{colortbl}
\usepackage{caption}
\usepackage[draft]{varioref}%
\providecommand{\U}[1]{\protect\rule{.1in}{.1in}}
\newtheorem{theorem}{Theorem}

\newtheorem{lemma}[theorem]{Lemma}

\newtheorem{remark}[theorem]{Remark}

\numberwithin{equation}{section}
\ifx\pdfoutput\relax\let\pdfoutput=\undefined\fi
\newcount\msipdfoutput
\ifx\pdfoutput\undefined\else
\ifcase\pdfoutput\else
\ifx\paperwidth\undefined\else
\ifdim\paperheight=0pt\relax\else\pdfpageheight\paperheight\fi
\ifdim\paperwidth=0pt\relax\else\pdfpagewidth\paperwidth\fi
\fi\fi\fi
\begin{document}

\title{The spectrum, radiation conditions and the Fredholm property for the Dirichlet
Laplacian in a perforated plane with semi-infinite inclusions}
\author{G. Cardone\thanks{University of Sannio, Department of Engineering, Corso
Garibaldi, 107, 82100 Benevento, Italy; email: giuseppe.cardone@unisannio.it},
T. Durante\thanks{University of Salerno, Department of Information
Engineering, Electrical Engineering and Applied Mathematics, Via Giovanni
Paolo II, 132, 84084 Fisciano (SA), Italy; email: tdurante@unisa.it.}, S.A.
Nazarov\thanks{St. Petersburg State University, 198504, Universitetsky pr.,
28, Stary Peterhof, Russia; Peter the Great St. Petersburg State Polytechnical
University, Polytechnicheskaya ul., 29, St. Petersburg, 195251, Russia;
Institute of Problems of Mechanical Engineering RAS, V.O., Bolshoj pr., 61,
St. Petersburg, 199178, Russia; email: srgnazarov@yahoo.co.uk.}}
\maketitle

\begin{abstract}
We consider the spectral Dirichlet problem for the Laplace operator in the
plane $\Omega^{\circ}$ with double-periodic perforation but also in the domain
$\Omega^{\bullet}$ with a semi-infinite foreign inclusion so that the
Floquet-Bloch technique and the Gelfand transform do not apply directly. We
describe waves which are localized near the inclusion and propagate along it.
We give a formulation of the problem with radiation conditions that provides a
Fredholm operator of index zero. The main conclusion concerns the spectra
$\sigma^{\circ}$ and $\sigma^{\bullet}$ of the problems in $\Omega^{\circ}$
and $\Omega^{\bullet},$ namely we present a concrete geometry which supports
the relation $\sigma^{\circ}\varsubsetneqq\sigma^{\bullet}$ due to a new
non-empty spectral band caused by the semi-infinite inclusion called an open
waveguide in the double-periodic medium.


\medskip


\end{abstract}

\section{Introduction\label{sect1}}

\subsection{Formulation of problems.\label{sect1.1}}

\begin{figure}
\begin{center}
\includegraphics[scale=0.45]{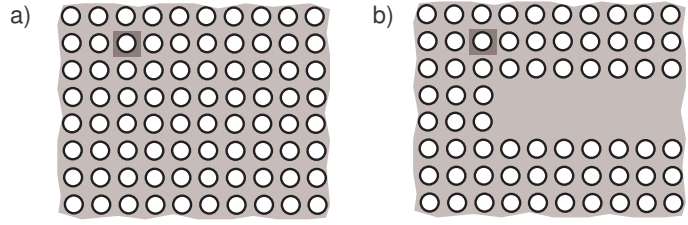}
\end{center}
\caption{The double-periodic perforated plane (a) and the semi-infinite
inclusion in it (b). The periodicity cell is shaded.}
\label{fig.1}
\end{figure}

Let $\omega\subset\mathbb{R}^{2}$ be a non-empty open set with smooth boundary
$\partial\omega$ such that the closure $\overline{\omega}=\omega\cup
\partial\omega$ belongs to the rectangle%
\begin{equation}%
\mathbb{Q}
=\left\{  x=\left(  x_{1},x_{2}\right)  :\left\vert x_{j}\right\vert
<l_{j},\ j=1,2\right\}  ,\ l_{j}>0.\label{1}%
\end{equation}
An infinite domain $\Omega^{\circ}$, fig. \ref{fig.1}, a, is the perforated
plane%
\begin{equation}
\Omega^{\circ}=\mathbb{R}^{2}\setminus%
{\textstyle\bigcup\limits_{\alpha\in\mathbb{Z}^{2}}}
\overline{\omega\left(  \alpha\right)  }\label{2}%
\end{equation}
where $\alpha=\left(  \alpha_{1},\alpha_{2}\right)  $, $\mathbb{Z}=\left\{
0,\pm1,\pm2,...\right\}  $ and%
\begin{equation}
\omega\left(  \alpha\right)  =\left\{  x:\left(  x_{1}-2\alpha_{1}l_{1}%
,x_{2}-2\alpha_{2}l_{2}\right)  \in\omega\right\}  .\label{3}%
\end{equation}%

Another domain $\Omega^{\bullet}$, fig. \ref{fig.1}, b, is obtained from
$\Omega^{\circ}$ by filling one $\left(  j=1\right)  $ or several
semi-infinite rows of holes (\ref{3}) with $\alpha_{1}\in\mathbb{N}=\left\{
1,2,3,...\right\}  $ and $\alpha_{2}=1,...,J$, that is,
\begin{equation}
\Omega^{\bullet}=\Omega^{\circ}\cup\Xi^{+}:=\Omega^{\circ}\cup\left\{
x_{1}:x_{1}>l_{1},x_{2}\in\left(  l_{2},2jl_{2}+l_{2}\right)  \right\}
.\label{4}%
\end{equation}
The spectral Dirichlet problem%
\begin{align}
-\Delta u\left(  x\right)   &  =\lambda u\left(  x\right)  ,\ x\in
\Omega^{\bullet},\label{5}\\
u\left(  x\right)   &  =0,\ x\in\partial\Omega^{\bullet},\label{6}%
\end{align}
and its weak formulation%
\begin{equation}
\left(  \nabla u,\nabla v\right)  _{\Omega^{\bullet}}=\lambda\left(
u,v\right)  _{\Omega^{\bullet}}\text{ \ }\forall v\in H_{0}^{1}\left(
\Omega^{\bullet}\right)  \label{5N}%
\end{equation}
are associated with an unbounded positive definite self-adjoint operator
$A^{\bullet}$ in $L^{2}\left(  \Omega^{\bullet}\right)  $ because the bilinear
form on the left-hand side of the integral identity (\ref{5N}) is a positive
definite form, closed in $H_{0}^{1}\left(  \Omega^{\bullet}\right)  ,$ see,
e.g., \cite[Ch 10]{BiSo}. Since the boundary $\partial\Omega^{\bullet}$ is
smooth, the domain of this operator becomes
\begin{equation}
\mathfrak{D}\left(  A^{\bullet}\right)  =H^{2}\left(  \Omega^{\bullet}\right)
\cap H_{0}^{1}\left(  \Omega^{\bullet}\right)  .\label{7}%
\end{equation}
In (\ref{5}) and (\ref{5N}), $\nabla=\operatorname{grad},$ $\Delta=\nabla
\cdot\nabla$ is the Laplace operator$,$ $\lambda$ a spectral parameter,
$(\cdot,\cdot)_{\Omega^{\bullet}}$ the natural scalar product in the Lebesgue
space $L^{2}\left(  \Omega^{\bullet}\right)  ,$ $H^{2}\left(  \Omega^{\bullet
}\right)  $ the Sobolev space, and $H_{0}^{1}\left(  \Omega^{\bullet}\right)
$ \ the subspace of functions $u\in H^{1}\left(  \Omega^{\bullet}\right)  $
satisfying the Dirichlet condition (\ref{6}).

The Dirichlet problem in the double-periodic domain (\ref{2}) is also supplied
with the operator $A^{\circ}$ in $L^{2}\left(  \Omega^{\circ}\right)  $
possessing the same general properties as $A^{\bullet}$. It is known, cf.
\cite{KuchUMN, Skrig, Kuchbook} and others, that the spectrum $\sigma^{\circ}$
of $A^{\circ}$ has the band-gap structure%
\begin{equation}
\sigma^{\circ}=%
{\textstyle\bigcup\limits_{n\in\mathbb{N}}}
B_{n}^{\circ} \label{8}%
\end{equation}
where the bands $B_{n}^{\circ}$, finite closed connected segments (\ref{16}),
will be described in Section \ref{sect2.1}.

According to \cite{CaNaTa}, the spectrum $\sigma^{\bullet}$ of the operator
$A^{\bullet}$ (and problems (\ref{5}), (\ref{6}) or (\ref{5N})) gets much more
complicated structure. One of goals in our paper is to find geometrical shapes
in (\ref{4}) such that the spectrum $\sigma^{\bullet}$ obtains at least one
additional band%
\begin{equation}
B_{0}^{\bullet}\subset\left(  0,\lambda_{\dag}^{\circ}\right]  \label{9}%
\end{equation}
below the cutoff point $\lambda_{\dag}^{\circ}=\underline{\sigma^{\circ}%
}:=\min\left\{  \lambda:\lambda\in\sigma^{\circ}\right\}  $ of the spectrum
(\ref{8}). However, the main purpose is to describe oscillatory waves which
are localized near the semi-infinite intact strip, cf. (\ref{4}),%
\begin{equation}
\Xi^{+}=\left(  l_{1},+\infty\right)  \times\left(  l_{2},\left(  2j+1\right)
l_{2}\right)  \label{10}%
\end{equation}
and travel along it. These waves decay exponentially as $x_{2}\rightarrow
\pm\infty$ and require for a radiation principle to detect direction of their
propagation, see Section \ref{sect5}. Moreover, radiation conditions provide
the problem (\ref{5}), (\ref{6}) with a Fredholm operator of index zero.

\begin{figure}
\begin{center}
\includegraphics[scale=0.45]{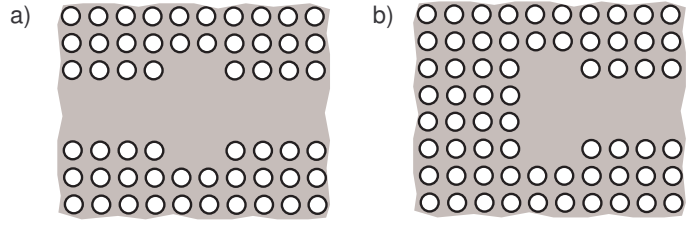}
\end{center}
\caption{Perturbed open waveguides supporting trapped modes.}
\label{fig.2}
\end{figure}

Using a primitive trick we also indicate geometries, fig. \ref{fig.2}, a and
b, which support trapped modes, i.e. eigenfunctions with the exponential decay
in all directions. Our approach can be readily adapted to other shapes of open
waveguides, see Section \ref{sect5.4}.%

\subsection{Preliminary discussion\label{sect1.2}}

\begin{figure}
\begin{center}
\includegraphics[scale=0.45]{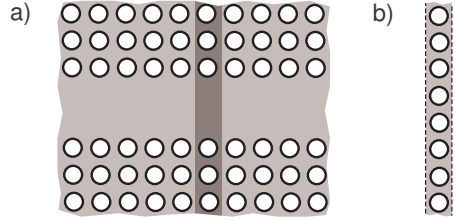}
\end{center}
\caption{The perforated plane with infinite inclusion (a) and the
corresponding infinite periodicity cell (b).}
\label{fig.3}
\end{figure}

Artificial experiments and natural phenomena demonstrate that waves may
propagate along rows of foreign inclusions in homogeneous and periodic
composite media. Classical mathematical tools to describe such wave processes
used to consider cases when corresponding boundary-value problems keep
periodicity at least in one direction, cf. the review papers \cite{BBDS, LM,
na417} and others. In our case the Dirichlet problem%
\begin{align}
-\Delta v\left(  x\right)  -\lambda v\left(  x\right)   &  =f\left(  x\right)
,\ x\in\Omega^{\sharp},\label{A1}\\
v\left(  x\right)   &  =0,\ x\in\partial\Omega^{\sharp}, \label{A2}%
\end{align}
must be posed in the domain, fig. \ref{fig.3}, a,%
\begin{equation}
\Omega^{\sharp}=\Omega^{\circ}\cup\Xi\label{A3}%
\end{equation}
with the infinite strip $\Xi=\mathbb{R}\times\left(  l_{2},\left(
2j+1\right)  l_{2}\right)  ,$ cf. (\ref{10}). Reducing size $l_{1}$ to $1/2$
by rescaling, we express the remaining periodicity along the $x_{1}$-axis as
follows:%
\begin{equation}
\Omega^{\sharp}=\left\{  x:\left(  x_{1}\pm1,x_{2}\right)  \in\Omega^{\sharp
}\right\}  . \label{A4}%
\end{equation}

Taking (\ref{A4}) into account, we apply the partial Gelfand transform
\cite{Gel}, see also \cite{na17} and \cite[\S 3.4]{NaPl},%
\begin{equation}
v\left(  x\right)  \mapsto V\left(  x;\zeta\right)  =\frac{1}{\sqrt{2\pi}%
}\underset{\alpha_{1}\in\mathbb{Z}}{\sum}e^{-i\zeta\alpha_{1}}v\left(
x_{1}+\alpha_{1},x_{2}\right)  \label{A5}%
\end{equation}
where $\zeta\in\left[  -\pi,\pi\right]  $ is the dual variable or the Floquet
parameter, and reduce (\ref{A1}), (\ref{A2}) to the parameter-dependent
problem%
\begin{align}
-\Delta V\left(  x;\zeta\right)  -M\left(  \zeta\right)  V\left(
x;\zeta\right)   &  =F\left(  x;\zeta\right)  ,\ x\in\Pi^{\sharp},\label{A6}\\
V\left(  x;\zeta\right)   &  =0,\ x\in\Gamma^{\sharp}, \label{A7}%
\end{align}
with the quasi-periodicity conditions%
\begin{equation}
V\left(  \frac{1}{2},x_{2};\zeta\right)  =e^{i\zeta}V\left(  -\frac{1}%
{2},x_{2};\zeta\right)  ,\ \ \ \frac{\partial V}{\partial x_{1}}\left(
\frac{1}{2},x_{2};\zeta\right)  =e^{i\zeta}\frac{\partial V}{\partial x_{1}%
}\left(  -\frac{1}{2},x_{2};\zeta\right)  ,\ \ \ x_{2}\in\mathbb{R},
\label{A8}%
\end{equation}
in the perforated strip which is shaded in fig. \ref{fig.3}, a, and redrawn in
fig. \ref{fig.3}, b,%
\begin{equation}
\Pi^{\sharp}=\left\{  x\in\Omega^{\sharp}:\left\vert x_{1}\right\vert
<1/2\right\}  . \label{A9}%
\end{equation}
Here, $M\left(  \zeta\right)  $ is a new notation for the spectral parameter,
the conditions (\ref{A8}) are imposed on the lateral sides of the strip and
$\Gamma^{\sharp}$ stands for the interior part of the boundary $\partial
\Pi^{\sharp},$%
\begin{equation}
\Gamma^{\sharp}=\left\{  x\in\partial\Pi^{\sharp}:\left\vert x_{1}\right\vert
<1/2\right\}  =\underset{\alpha_{2}\in\mathbb{Z}\setminus\left\{
1,...,J\right\}  }{%
{\textstyle\bigcup}
}\partial\omega\left(  0,\alpha_{2}\right)  . \label{A10}%
\end{equation}

Several fruitful approaches have been developed to indicate trapped modes,
namely solutions $V\in H^{2}\left(  \Pi^{\sharp}\right)  $ of the homogeneous
$\left(  F=0\right)  $ problem (\ref{A6})-(\ref{A8}), see again the review
papers \cite{BBDS, LM, na417} and many other publications. If $M\left(
\zeta\right)  $ is an eigenvalue of this problem in $\Pi^{\sharp}$ and
$V\left(  \cdot;\zeta\right)  $ is the corresponding eigenfunction, then the
Floquet wave%
\begin{equation}
v\left(  x;\zeta\right)  =e^{i\zeta\alpha_{1}}V\left(  x;\zeta\right)
,\ \ \left(  x_{1}-\alpha_{1},x_{2}\right)  \in\Pi^{\sharp},\ \alpha_{1}%
\in\mathbb{Z}, \label{A11}%
\end{equation}
becomes smooth in the $1$-periodic domain (\ref{A4}) due to the
quasi-periodicity conditions (\ref{A8}) and satisfies the homogeneous $\left(
f=0\right)  $ problem (\ref{A1}), (\ref{A2}). Moreover, it gains the
exponential decay as $x_2\rightarrow\pm\infty$ but oscillates in the $x_{1}%
$-direction. In other words, the wave is localized near the horizontal strip
$\Xi$ and propagates along it.

However, this direct and inherent way to detect localized propagative waves
does not work in a case when the periodicity in the $x_{1}$-coordinate is
disturbed even inside a finite volume, cf. fig. \ref{fig.2}, a, because the
partial Gelfand transform (\ref{A5}) no longer applies. On this issue, there
exists quite few results, e.g., \cite{BBDHa1, BBDHa2}, in particular about the
absence of trapped waves decaying in all directions. The paper \cite{CaNaTa}
provides a description of the essential spectrum of semi-infinite and broken
open waveguides but gives neither concrete examples of localized waves, nor
necessary radiation conditions but our paper partly eliminates these omissions.

In the sequel we display localized waves under certain restrictions, some of
which, especially shape and homogeneous structure of the open waveguide
(\ref{10}), can be easily avoided and have been introduced in order to
simplify demonstration. The principal requirement concerns the position of the
special spectral band $B_0^{\bullet}$ which gives rise to localized
waves, namely it is situated below the spectrum $\sigma^{\circ}$ in the
double-periodic domain (\ref{2}) while general results in \cite{CaNaTa} permit
for nucleation of new bands inside each of the spectral bands in
$\sigma^{\circ}.$ To provide the existence of the band (\ref{9}) we need the assumption (\ref{27}) 
below which means that the foreign semi-infinite inclusion filling some holes, see fig. \ref{fig.1}, b,
is sufficiently wide. This assumption imposes an upper bound for the spectral parameter
\begin{equation}
\lambda\in\left(  0,\lambda^{\sharp}\right)  \label{A0}%
\end{equation}
so that the method developed in this paper does not allow us to examine the whole spectrum $\sigma^\bullet$ but only its bottom part.

The key point in our analysis of problem (\ref{5}), (\ref{6}) in domain
(\ref{4}) in fig. \ref{fig.1}, b, is a localization weighted estimate derived
in Section \ref{sect5.1} that proves the exponential decay of a solution in
$\Omega^{\bullet}$ in each direction to infinity, except along the
semi-infinite inclusion $\Xi^{+}.$ A trick we use below to derive this
estimate, is based on integration by parts and simple algebraic operations but
works only in the case (\ref{A0}) where $\lambda^{\sharp}\in(0,\lambda_{\dag}^{\circ})$ is a certain bound, see
(\ref{27}), due to a technical reason. In this way, it remains an open
question to construct an elemental example of specific propagative waves with
the spectral parameter $\lambda$ inside a non-empty gap between the bands
$B_{n}^{\circ}$ and $B_{n+1}^{\circ}$ with $n\geq1$ in the spectrum (\ref{8}).
In \cite{NaTaRu} examples of arbitrarily many non-empty spectral bands for the
Dirichlet Laplacian in a double-periodic perforated plane are given.
Investigation of the spectral bands for other geometries of double-periodic two-dimensional
structures are performed in \cite{BaNa} and \cite{BaCaNaTa}.

\subsection{Structure of the paper\label{sect1.3}}

In Sections \ref{sect2} and \ref{sect3} we present some mainly known
information of the spectra of the Dirichlet problem in the domains
$\Omega^{\bullet}$, $\Omega^{\sharp}$ and Floquet waves localized near the
infinite inclusion $\Xi.$ In Section \ref{sect4} we derive asymptotics at
infinity of solutions to the inhomogeneous problem (\ref{A1}), (\ref{A2}).
Although we follow the scheme in \cite{na17}, \cite[\S \ 3.4]{NaPl}, we have
to repeat all arguments because the model problem (\ref{A6})-(\ref{A8}) is
posed in the infinite periodicity cell $\Pi^{\sharp}$, (\ref{A9}). We also
verify in Section \ref{sect2.3} that this problem with the spectral parameter
(\ref{A0}) supports just one trapped mode, i.e., an eigenfunction in
$H^{2}\left(  \Pi^{\sharp}\right)  $ with the exponential decay at infinity.

We start the last but central section with proving the localization estimate
which demonstrates that a solution of problem (\ref{A1}), (\ref{A2}) with the
permitted exponential growth in all directions in the plane actually decays in
all direction except along the semi-infinite inclusion $\Xi^{+}$. Together
with asymptotic formulas from Section \ref{sect4.4} which helps to detach the
above-mentioned localized Floquet waves, that estimate allows us to formulate
in Section \ref{sect5.2} radiation conditions which supply the operator of
problem in $\Omega^{\bullet}$ with index zero. It should be emphasized that
the classification "outgoing/ incoming" for waves propagating along the open
waveguide $\Xi^{+}$ is based on calculation of the Poynting vector and
application of the Mandelstam energy principle. We also construct in Section
\ref{sect5.3} the (right) parametrix for the operator of problem (\ref{A1}),
(\ref{A2}) in order to confirm its Fredholm property and in Section
\ref{sect5.4} we demonstrate the existence of trapped modes, that is,
eigenfunctions enjoying the exponential decay in all directions. We finish the
paper with mentioning available generalizations.

\section{Spectra\label{sect2}}

\subsection{The periodicity cell\label{sect2.1}}

In the framework of the Floquet-Block theory the Gelfand transform \cite{Gel},
see also \cite{KuchUMN, Skrig, Kuchbook} and others, applied to the Dirichlet
problem in $\Omega^{\circ}$, see (\ref{2}), generates the following spectral
problem in the perforated rectangle (\ref{1}), the periodicity cell $\varpi=%
\mathbb{Q}
\setminus\overline{\omega}$ shaded in fig. \ref{fig.1},%
\begin{align}
-\Delta U\left(  x;\eta\right)   &  =\Lambda\left(  \eta\right)  U\left(
x;\eta\right)  ,\ \ x\in\varpi,\label{11}\\
U\left(  x;\eta\right)   &  =0,\text{ \ }x\in\partial\omega, \label{12}%
\end{align}
with the quasi-periodicity conditions%
\begin{equation}
\left.  \frac{\partial^{p}U}{\partial x_{j}^{p}}\left(  x;\eta\right)
\right\vert _{x_{j}=l_{j}}=e^{i\eta_{j}l_{j}}\left.  \frac{\partial^{p}%
U}{\partial x_{j}^{p}}\left(  x;\eta\right)  \right\vert _{x_{j}=-l_{j}%
},\ \ p=0,1,\ \text{ }j=1,2. \label{13}%
\end{equation}
Here,%
\begin{equation}
\eta\in\left(  \eta_{1},\eta_{2}\right)  \in\mathbb{Y}=\left[  -\frac{\pi
}{2l_{1}},\frac{\pi}{2l_{1}}\right]  \times\left[  -\frac{\pi}{2l_{2}}%
,\frac{\pi}{2l_{2}}\right]  \label{14}%
\end{equation}
is the dual variable of the Gelfand transform, the Floquet parameter. Note
that we have set $l_{1}=1/2$ in Section \ref{sect1.2}. Problem (\ref{11}%
)-(\ref{13}) has the discrete spectrum composing the monotone unbounded
sequence%
\begin{equation}
0<\Lambda_{1}\left(  \eta\right)  \leq\Lambda_{2}\left(  \eta\right)
\leq...\leq\Lambda_{n}\left(  \eta\right)  \leq...\rightarrow+\infty\label{15}%
\end{equation}
where eigenvalues are listed according to their multiplicity. The functions
$\mathbb{Y}\ni\eta\mapsto\Lambda_{n}\left(  \eta\right)  $ are continuous and
$\pi l_{j}^{-1}$-periodic in $\eta_{j},$ $j=1,2,$ while the bands in (\ref{8})
are the connected, closed and finite segments%
\begin{equation}
B_{n}^{\circ}=\left\{  \Lambda_{n}\left(  \eta\right)  :\eta\in\mathbb{Y}%
\right\}  . \label{16}%
\end{equation}
It is known, see, e.g., \cite{KuchUMN, Skrig, Kuchbook}, that the union
$\sigma^{\circ}$ of the bands (\ref{16}) represents the whole spectrum of the
Dirichlet problem (\ref{5}), (\ref{6}) in the double-periodic domain
$\Omega^{\circ}.$

\subsection{The lower bound $\lambda_{\dag}^{\circ}$ of the spectrum
$\sigma^{\circ}$\label{sect2.2}}

The next assertion is a piece of the mathematical folklore and the authors do
not know the very origin of this result which is supplied with a condensed
proof for reader's convenience since it will be of further use in Lemmas
\ref{lemmaBBB} and \ref{lemmaE}. We, for example, refer to \cite{Vanninatan}
where a similar trick was used in homogenization.

\begin{lemma}
\label{lemmaA}There holds the relationship%
\begin{equation}
\Lambda_{1}\left(  \mathbf{0}\right)  <\Lambda_{1}\left(  \eta\right)  \text{
\ }\forall\eta\in\mathbb{Y},\ \ \eta\neq\mathbf{0}=(0,0). \label{17}%
\end{equation}

\end{lemma}

\textbf{Proof.} Since (\ref{13}) with $\eta=\mathbf{0}$ turns into the
periodicity conditions, by virtue of the strict maximum principle, the
principal eigenfunction $U_{1}\left(  x;\mathbf{0}\right)  $ can be fixed real
positive for $x\in\varpi$ as well as $\partial_{\nu}U_{1}\left(
x;\mathbf{0}\right)  <0$ for $x\in\partial\omega$ where $\partial_{\nu}$ is
the outward normal derivative on $\partial\varpi.$ Hence, $Z\left(
x;\eta\right)  =U_{1}\left(  x;\mathbf{0}\right)  ^{-1}U_{1}\left(
x;\eta\right)  $ is continuously differentiable in $\overline{\varpi}$ up to
the smooth boundary $\partial\omega$, in particular belongs to $H^{1}\left(
\varpi\right)  $. Then the integral identity serving for problem
(\ref{11})-(\ref{13}), assures that%
\begin{align}
\Lambda_{1}\left(  \eta\right) & \left\Vert U_{1}^{\eta};L^{2}\left(
\varpi\right)  \right\Vert ^{2}    =\left(  \nabla\left(  ZU_{1}^{0}\right)
,\nabla\left(  ZU_{1}^{0}\right)  \right)  _{\varpi}\label{18}\\
&  =\left(  Z\nabla U_{1}^{0},Z\nabla U_{1}^{0}\right)  _{\varpi}+\left(
Z\nabla U_{1}^{0},U_{1}^{0}\nabla Z\right)  _{\varpi}+\left(  U_{1}^{0}\nabla
Z,Z\nabla U_{1}^{0}\right)  _{\varpi}+\left(  U_{1}^{0}\nabla Z,U_{1}%
^{0}\nabla Z\right)  _{\varpi}\nonumber\\
&  =-\left(  Z\bigtriangleup U_{1}^{0},ZU_{1}^{0}\right)  _{\varpi}-2\left(
Z\nabla U_{1}^{0},U_{1}^{0}\nabla Z\right)  _{\varpi}\nonumber\\
&  +\left(  Z\nabla U_{1}^{0},U_{1}^{0}\nabla Z\right)  _{\varpi}+\left(
U_{1}^{0}\nabla Z,Z\nabla U_{1}^{0}\right)  _{\varpi}+\left\Vert U_{1}%
^{0}\nabla Z;L^{2}\left(  \varpi\right)  \right\Vert ^{2}\nonumber\\
&  =\Lambda_{1}\left(  \mathbf{0}\right)  \left\Vert U_{1}^{\eta};L^{2}\left(
\varpi\right)  \right\Vert ^{2}+\left\Vert U_{1}^{0}\nabla Z;L^{2}\left(
\varpi\right)  \right\Vert ^{2}+(\left(  U_{1}^{0}\nabla Z,Z\nabla U_{1}%
^{0}\right)  _{\varpi}-\left(  Z\nabla U_{1}^{0},U_{1}^{0}\nabla Z\right)
_{\varpi})\nonumber
\end{align}
where $U_{1}^{\eta}\left(  x\right)  =U_{1}\left(  x;\eta\right)  $. Being
pure imaginary, the latter difference vanishes because all other terms on the
left and right in (\ref{18}) are real. The function $Z=U_{1}^{\eta}/U_{1}^{0}$
cannot be constant in $\varpi$ for $\eta\neq\mathbf{0}$ due to the periodicity
of $U_{1}^{0}$ and the authentic quasi-periodicity of $U_{1}^{\eta}$. This concludes
with (\ref{17}). $\boxtimes$

We further need the principal eigenvalue $\Lambda^{\ast}>0$ of the mixed
boundary-value problem%
\begin{align}
-\Delta U\left(  x\right)   &  =\Lambda^{\star}U\left(  x\right)
,\ x\in\varpi,\ \ \ \ \partial_{\nu}U\left(  x\right)  =0,\ x\in\partial%
\mathbb{Q}
=\partial\varpi\setminus\partial\omega,\label{19}\\
U\left(  x\right)   &  =0,\ x\in\partial\omega, \label{20}%
\end{align}
together with the Friedrichs inequality%
\begin{equation}
\left\Vert \nabla U;L^{2}\left(  \varpi\right)  \right\Vert ^{2}\geq
\Lambda^{\star}\left\Vert U;L^{2}\left(  \varpi\right)  \right\Vert ^{2}\text{
\ }\forall U\in H_{0}^{1}\left(  \varpi,\partial\omega\right)  \label{21}%
\end{equation}
where $H_{0}^{1}\left(  \varpi,\partial\omega\right)  $ consists of functions
in $H^{1}\left(  \varpi\right)  $ which verify (\ref{20}). By the min
principle, cf. \cite[Thm. 10.2.1]{BiSo}, we have%
\begin{equation}
\lambda_{\dag}^{\circ}=\Lambda_{1}\left(  \mathbf{0}\right)  \geq
\Lambda^{\star}. \label{22}%
\end{equation}
Notice that $\Lambda_{1}\left(  \mathbf{0}\right)  =\Lambda^{\star}$ when the
cell $\varpi$ is symmetric with respect to both axes $x_{1}$ and $x_{2.}$

\subsection{Trapped modes\label{sect2.3}}

We fix some $\zeta\in\left[  -\pi,\pi\right]  $ and consider the Helmholtz
equation%
\begin{equation}
-\Delta V\left(  x;\zeta\right)  =M\left(  \zeta\right)  V\left(
x;\zeta\right)  ,\ x\in\Pi^{\sharp}, \label{23}%
\end{equation}
with the Dirichlet (\ref{A7}) and quasi-periodicity (\ref{A8}) conditions.
This problem is associated, see \cite[\S 10.1]{BiSo}, with a positive definite
self-adjoint operator $A^{\sharp}\left(  \zeta\right)  $ in $L^{2}\left(
\Pi^{\sharp}\right)  $ with the domain%
\begin{equation}
\mathfrak{D(}A^{\sharp}\left(  \zeta\right)  )=\left\{  V\in H^{2}\left(
\varpi\right)  :\text{(\ref{A7}) and the first equation in (\ref{A8}) are
met}\right\}  . \label{24}%
\end{equation}
According to \cite{na17}, see also \cite[\S 3.4]{NaPl}, the essential spectrum
$\sigma_{e}^{\sharp}\left(  \zeta\right)  $ of $A^{\sharp}\left(
\zeta\right)  $ takes the form%
\begin{equation}
\sigma_{e}^{\sharp}\left(  \zeta\right)  =\underset{n\in\mathbb{N}}{%
{\textstyle\bigcup}
}B_{n}^{\sharp}\left(  \zeta\right)  ,\text{ \ }B_{n}^{\sharp}\left(
\zeta\right)  =\{\Lambda_{n}\left(  \zeta,\eta_{2}\right)  :\eta_{2}\in
\lbrack-\tfrac{\pi}{l_{2}},\tfrac{\pi}{l_{2}}]\}. \label{25}%
\end{equation}
Hence, we recall that $l_{1}=1/2$ and, in view of (\ref{8}) and (\ref{16}),
write%
\begin{equation}
\sigma_{e}^{\circ}=\sigma^{\circ}=\underset{_{\zeta\in\left[  -\pi,\pi\right]
}}{%
{\textstyle\bigcup}
}\sigma_{e}^{\sharp}\left(  \zeta\right)  ,\text{ \ }B_{n}^{\circ}%
=\underset{\zeta\in\left[  -\pi,\pi\right]  }{%
{\textstyle\bigcup}
}B_{n}^{\sharp}\left(  \zeta\right)  . \label{26}%
\end{equation}
Using a standard argument, see \cite{Jones} and, e.g., \cite{na459}, we
examine the discrete spectrum $\sigma_{d}^{\sharp}\left(  \zeta\right)  $ of
$A^{\sharp}\left(  \zeta\right)  $ at $\zeta=0$ inside interval (\ref{A0}) and
then discuss the case $\zeta\neq0.$ In what follows we fix width $2l_{2}J$ of
the strip (\ref{10}) sufficiently large.

\begin{lemma}
\label{lemmaB}Under the condition%
\begin{equation}
\lambda^{\sharp}:=M^{\sharp}:=\pi^{2}l_{2}^{-2}\left(  2J\right)  ^{-2}%
<\min\left\{  \pi^{2},\Lambda^{\star}\right\}  \label{27}%
\end{equation}
the interval $(0,M^{\sharp})$ contains a unique eigenvalue $M_{1}\left(
0\right)  $, see (\ref{10}) and (\ref{22}), of the operator $A^{\sharp}\left(
0\right)  .$
\end{lemma}

\textbf{Proof.} We employ the max-min principle, cf. \cite[Thm 10.2.2]{BiSo},%
\begin{equation}
M_{p}\left(  \zeta\right)  =\underset{E_{p}^{\sharp}\left(  \zeta\right)
}{\max}\ \underset{V\in E_{p}^{^{\sharp}}\left(  \zeta\right)  \setminus
\left\{  0\right\}  }{\inf}\frac{\left\Vert \nabla V;L^{2}\left(  \Pi^{\sharp
}\right)  \right\Vert ^{2}}{\left\Vert V;L^{2}\left(  \Pi^{\sharp}\right)
\right\Vert ^{2}}, \label{28}%
\end{equation}
where $E_{p}^{\sharp}\left(  \zeta\right)  $ is any subspace of codimension
$p-1$ in the space
\begin{equation}
E^{\sharp}\left(  \zeta\right)  =\left\{  V\in H_{0}^{1}\left(  \Pi^{\sharp
};\Gamma^{\sharp}\right)  :V\left(  \tfrac{1}{2},x_{2}\right)  =e^{i\zeta
}V\left(  -\tfrac{1}{2},x_{2}\right)  ,\ x_{2}\in\mathbb{R}\right\}  .
\label{29}%
\end{equation}
Namely, if the right-hand side of (\ref{28}) with a $p\in\mathbb{N}$ is
strictly smaller than the lower bound $\sigma_{e}^{\sharp}\left(
\zeta\right)  $ of the essential spectrum in (\ref{25}), then the discrete
spectrum $\sigma_{d}^{\sharp}\left(  \zeta\right)  $ contains eigenvalues
$M_{1}\left(  \zeta\right)  ,...,M_{p}\left(  \zeta\right)  $ computed by
(\ref{28}). It should be emphasized that the space (\ref{29}) differs from
(\ref{24}) and, according to \cite[Ch.10]{BiSo}, coincides with the domain of
the bi-linear form $\left(  \nabla U,\nabla V\right)  _{\Pi^{\sharp}}$ in the
weak formulation of problem (\ref{A6})-(\ref{A8}). Taking $\zeta=0$ and $p=1$,
we insert into the Rayleigh quotient on the right-hand side of (\ref{28}) the
function%
\begin{equation}
V_{0}\left(  x\right)  =\sin\left(  \pi\left(  2J\right)  ^{-1}\left(
l_{2}^{-1}x_{2}-1\right)  \right)  \label{29N}%
\end{equation}
extended as null from $%
\mathbb{Q}
^{\sharp}=\Pi^{\sharp}\cap\Xi$ onto $\Pi^{\sharp}$. This function lives in
$E^{\sharp}\left(  0\right)  $ and makes the quotient equal to $M^{\sharp}$
from (\ref{26}). As a result, there exists an eigenvalue, $M_{1}\left(
0\right)  <M^{\sharp}<\Lambda^{\star}\leq\underline{\sigma_{e}^{\sharp}\left(
0\right)  }$, the lower bound of the essential spectrum, cf. (\ref{17}),
(\ref{22}) and (\ref{25}), while the first inequality is strict because in
contrast to our test function, an eigenfunction cannot vanish at a set of
positive area. Hence $M_{1}\left(  0\right)  \in\sigma_{d}^{\sharp}\left(
0\right)  .$

Let $p=2$ in (\ref{28}). The second, that is, first positive eigenvalue $N$ of
the Neumann problem in the rectangle $Q^{\sharp}=\left(  -1/2,1/2\right)
\times\left(  l_{2},\left(  2J+1\right)  l_{2}\right)  ,$ satisfies $N=\pi
^{2}\min\{1,\left(  2Jl_{2}\right)  ^{-2}\}=\pi^{2}\left(  2Jl_{2}\right)
^{-2},$ see (\ref{26}), while the ortogonality condition%
\begin{equation}
\int_{Q^{\sharp}}V\left(  x\right)  dx=0 \label{30}%
\end{equation}
assures the Poincar\'{e} inequality%
\begin{equation}
\int_{Q^{\sharp}}\left\vert \nabla V\left(  x\right)  \right\vert ^{2}dx\geq
N\int_{Q^{\sharp}}\left\vert V\left(  x\right)  \right\vert ^{2}dx. \label{31}%
\end{equation}
Adding to (\ref{31}) the Friedrichs inequalities (\ref{21}) in the cells
$\varpi\left(  0,\alpha_{2}\right)  \subset\Pi^{\sharp}$ with $\alpha_{2}%
\in\mathbb{Z}\setminus\left\{  1,...,J\right\}  ,$ we conclude that any
function $V$ in the subspace%
\begin{equation}
E^{\perp}\left(  0\right)  =\{V\in E^{\sharp}\left(  0\right)
:\text{(\ref{30}) is fulfilled}\} \label{32}%
\end{equation}
of codimension $1$ due to one orthogonality condition imposed, verifies the
estimate%
\begin{equation}
\left\Vert \nabla V;L^{2}\left(  \Pi^{\sharp}\right)  \right\Vert ^{2}\geq
\min\left\{  N,\Lambda^{\star}\right\}  \left\Vert V;L^{2}\left(  \Pi^{\sharp
}\right)  \right\Vert ^{2}. \label{33}%
\end{equation}
Thus, relations (\ref{27}), (\ref{33}) and (\ref{28}), $p=2$, show that the
second eigenvalue $M_{2}\left(  0\right)  $, if exists, lays outside the
interval $\left(  0,M^{\sharp}\right)  .\ \boxtimes$

In the case $\zeta\in\left(  0,\pi\right]  $ we replace the test function
(\ref{29N}) by $V_{\zeta}\left(  x\right)  =e^{i\zeta x_{1}}V_{0}\left(
x\right)  $ which clearly satisfies the first quasi-periodicity condition in
(\ref{A8}), compare with (\ref{29}). We have $\left\Vert V_{\zeta}%
;L^{2}\left(  \Pi^{\sharp}\right)  \right\Vert =\left\Vert V_{0};L^{2}\left(
\Pi^{\sharp}\right)  \right\Vert $ and%
\[
\left\Vert \nabla V_{\zeta};L^{2}\left(  \Pi^{\sharp}\right)  \right\Vert
^{2}=\left\Vert \nabla V_{0};L^{2}\left(  \Pi^{\sharp}\right)  \right\Vert
^{2}+\zeta^{2}\left\Vert V_{0};L^{2}\left(  \Pi^{\sharp}\right)  \right\Vert
^{2}.
\]
The max-min principle (\ref{28}), $p=1$, ensures the existence of an
eigenvalue $M_{1}\left(  \zeta\right)  \in\sigma_{d}^{\sharp}\left(
\zeta\right)  $ together with the estimate%
\[
M_{1}\left(  \zeta\right)  <M^{\sharp}+\zeta^{2}%
\]
but this inference is surely true only under the restriction $M^{\sharp}%
+\zeta^{2}<\underline{\sigma_{e}^{\sharp}\left(  \zeta\right)  }$. In view of
(\ref{25}) and Lemma \ref{lemmaA} the latter is valid for $\left\vert
\zeta\right\vert <\zeta^{\sharp}$ with some $\zeta^{\sharp}>0.$ Since
inequalities (\ref{31}) and (\ref{34}) do not take into account boundary
conditions at the lateral sides of the perforated strip $\Pi^{\sharp},$
formula (\ref{33}) does not involve the parameter $\zeta$ and, therefore, the
interval $\left(  0,M^{\sharp}\right)  $ may include at most one eigenvalue.

\begin{lemma}
\label{lemmaBBB}There exists $\zeta^{\sharp}>0$ such that, for $\left\vert
\zeta\right\vert <\zeta^{\sharp}$, the interval $\left(  0,M^{\sharp}\right)
$ contains the only eigenvalue $M_{1}\left(  \zeta\right)  $ in the discrete
spectrum $\sigma_{d}^{\sharp}\left(  \zeta\right)  $ of the operator
$A^{\sharp}\left(  \zeta\right)  .$ Moreover,%
\begin{equation}
M_{1}\left(  \zeta\right)  >M_{1}\left(  0\right)  \text{ \ for }\zeta
\in\left(  -\zeta^{\sharp},0\right)  \cup\left(  0,\zeta^{\sharp}\right).
\label{34N}%
\end{equation}

\end{lemma}

\textbf{Proof.} It remains to verify (\ref{34N}) and we apply the same
argument as in Lemma \ref{lemmaA}. Let $V_{1}^{\zeta}\left(  x\right)
=V_{1}\left(  x;\zeta\right)  $ be the eigenfunction of problem (\ref{23}),
(\ref{A7}), (\ref{A8}) corresponding to $M_{1}(\zeta).$ Again, by the strong
maximum principle, $V_{1}^{0}$ can be fixed positive in $\Pi^{\sharp}$ with
the negative outward normal derivative $\partial_{\nu}V_{1}^{0}$ on
$\Gamma^{\sharp}$. The fraction $Z^{\zeta}\left(  x\right)  =V_{1}^{0}\left(
x\right)  ^{-1}V_{1}^{\zeta}\left(  x\right)  $ is continuously differentiable
in $\overline{\Pi}^{\sharp}$ but does not belong to $H^{1}\left(  \Pi^{\sharp
}\right)  $. However, all integrals in the modified calculation (\ref{18}),%
\begin{align}
M_{1}(\zeta)||V_{1}^{\zeta};L^{2}\left(  \Pi^{\sharp}\right)  ||^{2}  &
=M_{1}(0)||V_{1}^{\zeta};L^{2}\left(  \Pi^{\sharp}\right)  ||^{2}+||V_{1}%
^{0}\nabla Z^{\zeta};L^{2}\left(  \Pi^{\sharp}\right)  ||^{2}\label{35N}\\
&  +\left(  \left(  V_{1}^{0}\nabla Z^{\zeta},Z^{\zeta}\nabla V_{1}%
^{0}\right)  _{\Pi^{\sharp}}-\left(  Z^{\zeta}\nabla V_{1}^{0},V_{1}^{0}\nabla
Z^{\zeta}\right)  _{\Pi^{\sharp}}\right) \nonumber\\
&  =M_{1}(0)||V_{1}^{\zeta};L^{2}\left(  \Pi^{\sharp}\right)  ||^{2}%
+||V_{1}^{0}\nabla Z^{\zeta};L^{2}\left(  \Pi^{\sharp}\right)  ||^{2}%
,\nonumber
\end{align}
converge owing to the exponential decay of $V_{1}^{0}\left(  x\right)  $ and
$V_{1}^{\zeta}\left(  x\right)  $ as $x_{2}\rightarrow\pm\infty,$ cf. Section
\ref{sect4.2}. The last norm in (\ref{35N}) is positive because $Z^{\zeta}$
cannot be constant in view of the periodicity of $V_{1}^{0}$ and the
quasi-periodicity of $V_{1}^{\zeta}.\ \boxtimes$

\bigskip

\begin{figure}
\begin{center}
\includegraphics[scale=0.45]{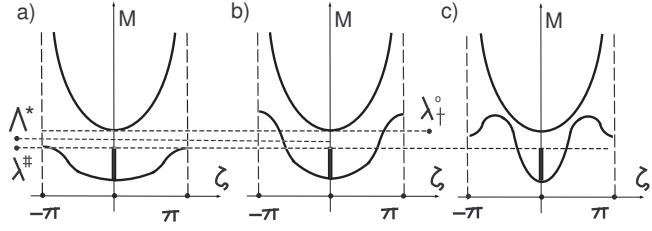}
\end{center}
\caption{Several possible positions of the curve $M=M_{1}(\zeta)$.}
\label{fig.4}
\end{figure}

Two typical dispositions of the eigenvalue $M_{1}(\zeta)$ below the essential
spectrum are depicted in fig. \ref{fig.4}, a and b, while $\zeta^{\sharp}=\pi$
in the first case but $\zeta^{\sharp}<\pi$ in the second one.%

When $M_{1}(\zeta)$ stays in the discrete spectrum, the function $\zeta\mapsto
M_{1}(\zeta)$ is continuous and even (the latter is verified by complex
conjugation in problem (\ref{23}), (\ref{A7}), (\ref{A8}) so that the graphs
in fig. \ref{fig.4} are symmetric with respect to the ordinate axis). Hence,
the set
\begin{equation}
B^{\sharp}=\left\{  M_{1}(\zeta):\zeta\in\left(  -\zeta_{0},\zeta_{0}\right)
\right\}  =[M_{1}(0),M^{\sharp}) \label{MMM}%
\end{equation}
with $M^{\sharp}=\pi^{2}(2l_{2}J)^{-1}$, (\ref{27}), is a semi-open segment.
According to \cite{CaNaTa}, this set is just a part of the lowest additional
segment (\ref{9}) in the essential spectrum $\sigma_{e}^{\bullet}$ of problem
(\ref{5}), (\ref{6}).

We are not able to reject the graph $M=M_{1}(\zeta)$ with several local
extrema, cf. fig. \ref{fig.4}, c, but the upper bound $M^{\sharp}$ in
(\ref{MMM}) is fixed such that, for any $M\in(0,M^{\sharp})\subset
B_{0}^{\bullet}$ (the thick line in fig. \ref{fig.4}, a-c), there exist
exactly two points $\pm\zeta\left(  M\right)  $ with $M\left(  \pm\zeta\left(
M\right)  \right)  =M$. This restriction is introduced to simplify further
notation in Section \ref{sect5}.

\section{Localized Floquet waves\label{sect3}}

\subsection{Propagative, standing and resonance waves\label{sect3.1}}

First of all, we make the change%
\begin{equation}
V\left(  x;\zeta\right)  \mapsto W\left(  x;\zeta\right)  =e^{-i\zeta x_{1}%
}V\left(  x;\zeta\right)  \label{34}%
\end{equation}
and rewrite problem (\ref{23}), (\ref{A7}), (\ref{A8}) as follows:%
\begin{gather}
-\left(  \frac{\partial}{\partial x_{1}}+i\zeta\right)  ^{2}W\left(
x;\zeta\right)  -\frac{\partial^{2}}{\partial x_{2}^{2}}W\left(
x;\zeta\right)  =M(\zeta)W\left(  x;\zeta\right)  ,\ \ x\in\Pi^{\sharp
},\label{35}\\
W\left(  x;\zeta\right)  =0,\ x\in\Gamma^{\sharp},\label{36}\\
W\left(  \frac{1}{2},x_{2};\zeta\right)  =W\left(  -\frac{1}{2},x_{2}%
;\zeta\right)  ,\text{ }\frac{\partial W}{\partial x_{1}}\left(  \frac{1}%
{2},x_{2};\zeta\right)  =\frac{\partial W}{\partial x_{1}}\left(  -\frac{1}%
{2},x_{2};\zeta\right)  ,\ \ x\in\mathbb{R}\text{.}\label{37}%
\end{gather}
Now the boundary-value problem (\ref{35})-(\ref{37}) can be interpreted as a
polynomial, actually quadratic, pencil in the complex variable $\zeta
\in\mathbb{C}$, that is,
\begin{equation}
\mathbb{C}\ni\zeta\mapsto\left(  \mathfrak{A}^{\sharp}\left(  \zeta;M\right)
:H_{per}^{1}\left(  \Pi^{\sharp}\right)  \cap H_{0}^{1}\left(  \Pi^{\sharp
};\Gamma^{\sharp}\right)  \rightarrow L^{2}\left(  \Pi^{\sharp}\right)
\right)  ,\label{39}%
\end{equation}
see \cite{GoKr} and, e.g., \cite[Ch.1]{NaPl} for summary of results. Notice
that change (\ref{34}) purposes to make the domain of the pencil independent
of $\zeta.$

If, for some $\zeta\in\lbrack-\pi,\pi],$ $M$ is an eigenvalue of problem
(\ref{23}), (\ref{A7}), (\ref{A8}) and, therefore, of problem (\ref{35}%
)-(\ref{37}), then $\zeta$ is an eigenvalue of the pencil $\mathfrak{A}%
^{\sharp}\left(  \cdot;M\right)  $ with the same eigenfunction, namely an
eigenvector of $\mathfrak{A}\left(  \zeta;M\right)  $. However, the pencil may
get associated vectors $W^{1},...,W^{\varkappa-1}$ in addition to an
eigenvector $W^{0}$ which all together form a Jordan chain and have to be
found out from the abstract equations%
\begin{equation}
\mathfrak{A}^{\sharp}\left(  \zeta;M\right)  W^{k}=-\overset{k}{\underset
{l=1}{\sum}}\frac{1}{l!}\frac{\partial^{l}\mathfrak{A}^{\sharp}}{\partial
\zeta^{l}}\left(  \zeta;M\right)  W^{k-l},\text{ \ }k=0,...,\varkappa-1.
\label{40}%
\end{equation}
Taking $\zeta=0$ and $M=M_{1}\left(  0\right)  ,$ we write the corresponding
boundary-value problem at $k=1$ as the differential equation%
\begin{equation}
-\Delta W^{1}\left(  x\right)  -M_{1}\left(  0\right)  W^{1}\left(  x\right)
=F^{1}\left(  x\right)  :=2i\frac{\partial W^{0}}{\partial x_{1}}\left(
x\right)  ,x\in\Pi^{\sharp}, \label{41}%
\end{equation}
with the Dirichlet (\ref{36}) and the periodicity (\ref{37}) conditions. Since
$M_{1}\left(  0\right)  $ is a simple eigenvalue due to Lemma \ref{lemmaBBB}
and the problem is formally self-adjoint, the Fredholm alternative brings the
only compatibility condition in problem (\ref{41}), (\ref{36}), (\ref{37})%
\begin{equation}
\int_{\Pi^{\sharp}}W^{0}\left(  x\right)  F^{1}\left(  x\right)  dx=0
\label{400}%
\end{equation}
which is easily verified by integration by parts because $W^{0}$ is real.
Thus, the problem admits a solution which is determined up to an addendum
$cW^{0}$ and can be made pure imaginary. The constructed Jordan chain
$\left\{  W^{0},W^{1}\right\}  $ gives rise to two Floquet waves, standing and
resonance, namely bounded and with the linear growth,
\begin{align}
w^{0}\left(  x\right)   &  =V_{1}\left(  x;0\right)  =W^{0}\left(  x\right)
,\label{42}\\
w^{1}\left(  x\right)   &  =ix_{1}W^{0}\left(  x\right)  +W^{1}\left(
x\right)  \label{42N}%
\end{align}
which satisfy the homogeneous $\left(  f=0\right)  $ problem (\ref{A1}),
(\ref{A2}) in the periodic domain (\ref{A3}), fig. \ref{fig.3}, a.

The associated vector $W^{2}$ of rank $2$ must fulfil the differential
equation%
\begin{equation}
-\Delta W^{2}\left(  x\right)  -M_{1}\left(  0\right)  W^{2}\left(  x\right)
=F^{2}\left(  x\right)  =2i\frac{\partial W^{1}}{\partial x_{1}}\left(
x\right)  -W^{0}\left(  x\right)  ,\ \ x\in\Pi^{\sharp}, \label{43}%
\end{equation}
with the usual conditions (\ref{35}) and (\ref{36}). This is nothing but the
differential form of the abstract equation (\ref{40}) with $k=2$. Note that
the right-hand sides of (\ref{41}) and (\ref{43}) involve the first and
second-order derivatives of $\mathfrak{A}\left(  \zeta;M_{1}\left(  0\right)
\right)  $ in the variable $\zeta$. In the next lemma we will prove that this
problem has no solution and, therefore, the Jordan chain $\left\{  W^{0}%
,W^{1}\right\}  $ cannot be extended and there is no Floquet wave at
$M=M_{1}\left(  0\right)  $ with the quadratic growth as $x_{1}\rightarrow
\pm\infty$.

\begin{lemma}
\label{lemmaE}There holds the formula%
\begin{equation}
b:=-\left(  F^{2},W^{0}\right)  _{\Pi^{\sharp}}>0. \label{44}%
\end{equation}

\end{lemma}

\textbf{Proof.} We recall that $W^{0}\left(  x\right)  =V_{1}\left(
x,0\right)  >0,$ $x\in\Pi^{\sharp},$ and, hence, the function $Z=\left(
W^{0}\right)  ^{-1}W^{1}$ is pure imaginary and continuously differentiable
while all integrals below converge due to the exponential decay of $W^{0}$ and
$W^{1}$ as $x_{2}\rightarrow\pm\infty$, see Remark \ref{RemarkDECAY}. Equation
(\ref{41}) turns into%
\[
-\nabla\cdot\left(  W^{0}\nabla Z\right)  -\nabla W^{0}\cdot\nabla
Z=2i\frac{\partial W^{0}}{\partial x_{1}}\text{ \ in }\Pi^{\sharp}.
\]
Multiplying it with $W^{0}\overline{Z}$ and integrating by parts with the help
of the boundary conditions yield
\begin{align}
\left(  W^{0}\nabla Z,W^{0}\nabla Z\right)  _{\Pi^{\sharp}}  &  =2i\left(
\partial_{1}W^{0},W^{0}Z\right)  _{\Pi^{\sharp}}\label{Z11}\\
&  =i\left(  \partial_{1}W^{0},W^{0}Z\right)  _{\Pi^{\sharp}}-i\left(
W^{0},\partial_{1}\left(  W^{0}Z\right)  \right)  _{\Pi^{\sharp}}=-i\left(
W^{0},W^{0}\partial_{1}Z\right)  _{\Pi^{\sharp}},\nonumber
\end{align}
where $\partial_{j}=\partial/\partial x_{j}$. Furthermore,
\begin{align}
b  &  =\left(  W^{0},W^{0}\right)  _{\Pi^{\sharp}}-2i\left(  Z\partial
_{1}W^{0},W^{0}\right)  _{\Pi^{\sharp}}-2i\left(  W^{0}\partial_{1}%
Z,W^{0}\right)  _{\Pi^{\sharp}}\label{Z12}\\
&  =\left(  W^{0},W^{0}\right)  _{\Pi^{\sharp}}-i(Z,\partial_{1}\left(
W^{0}\right)  ^{2})_{\Pi^{\sharp}}-2i\left(  W^{0}\partial_{1}Z,W^{0}\right)
_{\Pi^{\sharp}}\nonumber\\
&  =\left(  W^{0},W^{0}\right)  _{\Pi^{\sharp}}+i(\partial_{1}Z,\left(
W^{0}\right)  ^{2})_{\Pi^{\sharp}}-2i\left(  W^{0}\partial_{1}Z,W^{0}\right)
_{\Pi^{\sharp}}.\nonumber
\end{align}
Combining (\ref{Z11}) and (\ref{Z12}), we have
\begin{align}
b  &  =\left(  W^{0},W^{0}\right)  _{\Pi^{\sharp}}-i\left(  W^{0}\partial
_{1}Z,W^{0}\right)  _{\Pi^{\sharp}}+i\left(  W^{0},W^{0}\partial_{1}Z\right)
_{\Pi^{\sharp}}+\left(  W^{0}\partial_{1}Z,W^{0}\partial_{1}Z\right)
_{\Pi^{\sharp}}\label{FFF}\\
&  =\left\Vert W^{0}\left(  \partial_{1}Z+i\right)  ;L^{2}\left(  \Pi^{\sharp
}\right)  \right\Vert ^{2}+\left\Vert W^{0}\partial_{2}Z;L^{2}\left(
\Pi^{\sharp}\right)  \right\Vert ^{2}>0.\nonumber
\end{align}
The strict inequality is valid because the $1-$periodic in $x_{1}$ function
$Z$ cannot be equal to $-ix_{1}$. $\boxtimes$

\bigskip

If $M=M_{1}\left(  \zeta\right)  \in B^{\sharp}$ and $0<\left\vert
\zeta\right\vert <\zeta^{\sharp}$, then problem (\ref{A1}), (\ref{A2}) has the
Floquet waves%
\begin{equation}
w^{\pm}\left(  x,\zeta\right)  =e^{\pm i\zeta x_{1}}W\left(  x\right)  ,
\label{44bis}%
\end{equation}
cf. (\ref{A11}) and (\ref{34}), where $W^{+}$ is an eigenfunction of problem
(\ref{35})-(\ref{37}) with $M=M_{1}\left(  \zeta\right)  $, $\zeta\in\left(
0,\pi\right)  $ and $W^{-}\left(  x\right)  =\overline{W^{+}\left(  x\right)
}.$ In contrast to the staying wave (\ref{42}), which is just $1$-periodic in
$x_{1}$, waves (\ref{44bis}) are oscillatory due to the factors $e^{\pm i\zeta
x_{1}}$ with the period $2\pi/\zeta\neq1$. The absence of other Floquet waves
occurs by virtue of our assumption in Section \ref{sect2.3} on the semi-open
segment (\ref{MMM}). A different argument will confirm this fact in Section
\ref{sect5.3}.

\begin{remark}
\label{RemarkONLY}The graph in fig. \ref{fig.4}, c, which is probably
impossible, furnishes four Floquet waves (\ref{44bis}) for $M\in\left(
M\left(  \pm\pi\right)  ,M^{\cap}\right)  $ where%
\begin{equation}
M^{\cap}=M\left(  \zeta^{\cap}\right)  =\max M\left(  \zeta\right)
,\ \ \zeta^{\cap}\in\left(  0,\pi\right)  . \label{cap}%
\end{equation}
Moreover, at $M=M^{\cap}$ there appear two linear Floquet waves
\begin{equation}
e^{\pm i\zeta^{\cap}x_{1}}\left(  \pm ix_{1}W_{\cap}^{0\pm}\left(  x\right)
+W_{\cap}^{1\pm}\left(  x\right)  \right)  \label{Lin}%
\end{equation}
in addition to waves (\ref{44bis}) with $\zeta=\zeta^{\cap}$ and $W^{\pm
}=W_{\cap}^{0\pm}$. To avoid the incipient inconsistency in the notation, we
had introduced a supplementary restriction on $M$ in Section \ref{sect2.3}
and, in particular consider the semi-open segment $B^{\sharp}$ even in the
case $\zeta^{\sharp}=\pi$ because the linear Floquet waves (\ref{Lin}) are
attributed to the point $\zeta=\pm\pi$ in fig. \ref{fig.4}, a, too. At the
same time, a clear modification of our notation is only needed to cover all
the discarded situations. $\boxtimes$
\end{remark}

\subsection{Outgoing and incoming waves\label{sect3.2}}

In \cite{na128, na147}, \cite[\S 5.4]{NaPl} it was shown that the\ symplectic,
that is, sesquilinear and anti-Hermitian form%
\begin{equation}
q_{R}\left(  u,v\right)  =\int_{\left\{  x\in\Omega^{\sharp}:x_{1}=R\right\}
}\left(  \overline{v\left(  R,x_{2}\right)  }\frac{\partial u}{\partial x_{1}%
}\left(  R,x_{2}\right)  -u\left(  R,x_{2}\right)  \overline{\frac{\partial
v}{\partial x_{1}}\left(  R,x_{2}\right)  }\right)  dx_{2} \label{F1}%
\end{equation}
is proportional to the mean-value of the projection on the $x_{1}$-axis of the
Poynting vector \cite{Poynt} which indicates the direction of energy transfer
by a propagative wave. We use this observation to classify Floquet waves
according to the Mandelstam energy radiation principle \cite{Mand}, see also
\cite[Ch.1]{VoBa}, \cite{na543}, \cite{na569} and others.

Since the form $q_{R}\left(  u,v\right)  $ appears as a line integral in the
Green formula for the Helmholtz operator and, hence, does not depend on the
parameter $R>0$ if $u$ and $v$ satisfy the homogeneous $\left(  f=0\right)  $
problem (\ref{A1}), (\ref{A2}). Moreover, after integration in $R\in\left(
n,n+1\right)  $ we obtain%
\begin{equation}
q\left(  u,v\right)  =\int_{\Pi^{\sharp}\left(  n\right)  }\left(
\overline{v\left(  x\right)  }\frac{\partial u}{\partial x_{1}}\left(
x\right)  -u\left(  x\right)  \overline{\frac{\partial v}{\partial x_{1}%
}\left(  x\right)  }\right)  dx, \label{F2}%
\end{equation}
where $\Pi^{\sharp}\left(  n\right)  =\{x:\left(  x_{1}-n,x_{2}\right)  \in
\Pi^{\sharp}\}$ is a shifted perforated strip (\ref{A9}) and $n\in\mathbb{N}.$

According to the Mandelstam energy principle \cite{Mand} in the interpretation
\cite{na128, na147}, \cite[\S 5.3]{NaPl}, we call a wave $w$ outgoing to
infinity in the case $\operatorname{Im}q\left(  w,w\right)  >0$ and incoming
from infinity in the case $\operatorname{Im}q\left(  w,w\right)  <0$.

For the Floquet waves (\ref{44bis}), we have%
\begin{align}
q\left(  w^{\pm},w^{\pm}\right)   &  =\int_{\Pi^{\sharp}}\left(
\overline{W^{\pm}\left(  x\right)  }\frac{\partial W^{\pm}}{\partial x_{1}%
}\left(  x\right)  \pm i\zeta W^{\pm}\left(  x\right)  -W^{\pm}\left(
x\right)  \left(  \overline{\frac{\partial W^{\pm}}{\partial x_{1}}\left(
x\right)  \pm i\zeta W^{\pm}\left(  x\right)  }\right)  \right)
dx\label{F3}\\
&  =2i\operatorname{Im}\int_{\Pi^{\sharp}}\overline{W^{\pm}\left(  x\right)
}\left(  \frac{\partial W^{\pm}}{\partial x_{1}}\left(  x\right)  \pm i\zeta
W^{\pm}\left(  x\right)  \right)  dx=ia^{\pm}.\nonumber
\end{align}
Recalling the formula $w^{+}=\overline{w^{-}},$ we see that $a^{\pm}=\pm a$
and $q\left(  w^{\pm},w^{\pm}\right)  =0$. It should be underlined that the
equality $a=0$ means that the differential equation
\begin{equation}
-\left(  \left(  \frac{\partial}{\partial x_{1}}\pm i\zeta\right)  ^{2}%
+\frac{\partial^{2}}{\partial x_{2}^{2}}+M\left(  \zeta\right)  \right)
W_{1}^{\pm}\left(  x\right)  =2i\left(  \frac{\partial}{\partial x_{1}}\pm
i\zeta\right)  W^{\pm}\left(  x\right)  ,\ \ x\in\Pi^{\sharp}, \label{F9}%
\end{equation}
with conditions (\ref{36}), (\ref{37}) is solvable and, hence, the eigenvalue
$\pm\zeta$ of the pencil $\mathfrak{A}\left(  \cdot;M\left(  \zeta\right)
\right)  $ has the associated vector $W_{1}^{\pm}$ in addition to the
eigenvector $W^{\pm}.$ For $\lambda\in(M_{1}\left(  0\right)  ,M^{\sharp})$,
the latter contradicts Lemma \ref{lemmaB} and our definition of the upper
bound $M^{\sharp}.$ In the next section we will show that $a^{+}>0$ if and
only if the function $\zeta\mapsto M\left(  \zeta\right)  $ is strictly
growing at the point $+\zeta$ so that in accord with fig. \ref{fig.4}, a and
b, the waves $w^{+}$ and $w^{-}$ are outgoing and incoming respectively.

The waves (\ref{42}) and (\ref{42N}) at $\lambda=M_{1}\left(  0\right)  $
satisfy $q\left(  w^{0},w^{0}\right)  =q\left(  w^{1},w^{1}\right)  =0$
because $w^{0}$ and $iw^{1}$ are fixed real. At the same time, we derive from
(\ref{44}), (\ref{43}) that%
\begin{align}
q\left(  w^{1},w^{0}\right)   &  =\int_{\Pi^{\sharp}}\left(  \overline
{W^{0}\left(  x\right)  }\left(  ix_{1}\frac{\partial W^{0}}{\partial x_{1}%
}\left(  x\right)  +\frac{\partial W^{1}}{\partial x_{1}}\left(  x\right)
+iW^{0}\left(  x\right)  \right)  \right. \label{F4}\\
&  -\left.  \left(  ix_{1}W^{0}\left(  x\right)  +W^{1}\left(  x\right)
\right)  \overline{\frac{\partial W^{0}}{\partial x_{1}}\left(  x\right)
}\right)  dx\nonumber\\
&  =\int_{\Pi^{\sharp}}\left(  W^{0}\left(  x\right)  \left(  \frac{\partial
W^{1}}{\partial x_{1}}\left(  x\right)  +iW^{0}\left(  x\right)  \right)
-W^{1}\left(  x\right)  \frac{\partial W^{0}}{\partial x_{1}}\left(  x\right)
\right)  dx=ib.\nonumber
\end{align}
Following \cite{na128, na147} and \cite[\S 5.3]{NaPl}, we introduce the linear
wave packets%
\begin{equation}
w^{\pm}\left(  x;0\right)  =w^{1}\left(  x\right)  \pm w^{0}\left(  x\right)
\label{66}%
\end{equation}
and recognize $w^{+}$ and $w^{-}$ are outgoing and incoming, respectively,
because $b>0$ in (\ref{F4}) due to Lemma \ref{lemmaE} and
\begin{equation}
q\left(  w^{\pm},w^{\pm}\right)  =\pm2ib,\text{ \ }q\left(  w^{\pm},w^{\mp
}\right)  =0. \label{F5}%
\end{equation}
We have classified the Floquet waves (\ref{66}) and (\ref{44bis}) for the
parameter $\lambda\in(M_{1}\left(  0\right)  ,M^{\sharp})$ as in fig.
\ref{fig.4}, a and b.

\subsection{Some asymptotic formulas\label{sect.3.3}}

Let
\begin{equation}
M^{\varepsilon}=M^{0}+\varepsilon\label{P1}%
\end{equation}
where $M^{0}\in B^{\sharp}$ and $\varepsilon>0$. We first set $M^{0}%
=M_{1}\left(  0\right)  $ and accept the standard asymptotic ans\"{a}tze
\cite[Ch. 9]{VaTr} for eigenvalues and eigenvectors of operator
pencils\footnote{The book \cite{VaTr} deals with linear but non self-adjoint
pencils and reduction of our quadratic pencil to that one is obvious.}%
\begin{align}
\zeta^{\varepsilon\pm}  &  =\pm\varepsilon^{1/2}\zeta^{\prime}+\varepsilon
\zeta^{\prime\prime}+...,\label{P2}\\
W^{\varepsilon\pm}\left(  x\right)   &  =W^{0}\left(  x\right)  \pm
\varepsilon^{1/2}\zeta^{\prime}W^{1}\left(  x\right)  +\varepsilon
W^{^{\prime\prime}}\left(  x\right)  +... \label{P3}%
\end{align}
where denominator in the exponent $1/2$ of the parameter $\varepsilon$ is
nothing that length $2$ of the Jordan chain. We insert them into the equation
\begin{equation}
\mathfrak{A}^{\sharp}(\zeta^{^{\varepsilon\pm}};M^{\varepsilon})W^{\varepsilon
\pm}=0 \label{P4}%
\end{equation}
and collect coefficients of $1,\pm\zeta^{\prime},\varepsilon^{1/2}$ and
$\varepsilon$. We obtain
\begin{align}
\mathfrak{A}^{\sharp}\left(  0;M^{0}\right)  W^{0}  &  =0,\text{
\ }\mathfrak{A}^{\sharp}\left(  0;M^{0}\right)  W^{1}=-\frac{d\mathfrak{A}%
^{\sharp}}{d\zeta}\left(  0;M^{0}\right)  W^{0},\label{P5}\\
\mathfrak{A}^{\sharp}\left(  0;M^{0}\right)  W^{^{\prime\prime}}  &
=-\frac{d\mathfrak{A}^{\sharp}}{d\zeta}\left(  0;M^{0}\right)  \left(
\zeta^{\prime\prime}W^{0}\pm\zeta^{\prime}W^{1}\right)  -\frac{1}{2}\left(
\zeta^{\prime}\right)  ^{2}\frac{d^{2}\mathfrak{A}^{\sharp}}{d\zeta^{2}%
}\left(  0;M^{0}\right)  W^{0}+W^{0} \label{P6}%
\end{align}
while the last term is due to the perturbation $\varepsilon$ in (\ref{P1}).

Equations (\ref{P5}) are nothing but (\ref{40}) with $k=0,1$ and, by
definition, elements of \ the Jordan chain fulfil them. In view of
(\ref{400}), (\ref{41}) and (\ref{44}), (\ref{43}) the compatibility condition
in equation (\ref{P6}) reads%
\[
-\left(  \zeta^{\prime}\right)  ^{2}b+\left\Vert W^{0};L^{2}\left(
\varpi\right)  \right\Vert ^{2}=0.
\]
The positive root of this quadratic equation%
\[
\zeta^{\prime}=b^{-1/2}\left\Vert W^{0};L^{2}\left(  \varpi\right)
\right\Vert
\]
specifies the main terms of ans\"{a}tze (\ref{P2}) and (\ref{P3}) while
general results in \cite[Ch.9]{VaTr} provide estimates of the asymptotic
remainders. The obtained formula%
\[
\left\vert \zeta^{\varepsilon\pm}\mp\varepsilon^{1/2}\zeta^{\prime}\right\vert
\leq c\varepsilon
\]
together with (\ref{P1}) imply that in the vicinity of the point $\left(
0,M_{1}\left(  0\right)  \right)  $ the graph of the function $M_{1}\left(
\zeta\right)  $ is approximated by the parabola $M_{1}\left(  0\right)
+b\left\Vert W^{0};L^{2}\left(  \varpi\right)  \right\Vert ^{-2}\zeta^{2}$ as
it is depicted in fig. \ref{fig.4}.

Let now $M^{0}\in(M_{1}\left(  0\right)  ,M^{\sharp})$. Then the pencil
$\mathfrak{A}^{\sharp}\left(  \cdot;M^{0}\right)  $ has two simple eigenvalues
$\pm\zeta^{0}$ with $\zeta^{0}\in\left(  0,\pi\right)  $ and the corresponding
eigenfunctions are denoted by $W^{\pm}$ while $W^{+}=\overline{W^{-}}.$ The
asymptotic ans\"{a}tze from the book \cite[Ch.9]{VaTr}%
\begin{equation}
\zeta^{\varepsilon\pm}=\pm\zeta^{0}\pm\varepsilon\zeta^{\prime}+...,\text{
\ }W^{\varepsilon\pm}\left(  x\right)  =W^{\pm}\left(  x\right)  +\varepsilon
W^{\pm\prime}\left(  x\right)  +... \label{P7}%
\end{equation}
inserted into (\ref{P4}) leads to the abstract equation%
\[
\mathfrak{A}^{\sharp}\left(  \pm\zeta^{0};M^{0}\right)  W^{\pm^{\prime}%
}=-\frac{d\mathfrak{A}^{\sharp}}{d\zeta}\left(  \pm\zeta^{0};M^{0}\right)
W^{\pm}+W^{\pm}%
\]
which turns into the differential equation%
\[
-\left(  \left(  \frac{\partial}{\partial x_{1}}\pm i\zeta^{0}\right)
^{2}+\frac{\partial^{2}}{\partial x_{2}^{2}}+M^{0}\right)  W^{\pm\prime
}\left(  x\right)  =\pm2i\zeta^{\prime}\left(  \frac{\partial}{\partial x_{1}%
}\pm i\zeta^{0}\right)  W^{\pm}\left(  x\right)  +W^{\pm}\left(  x\right)
,\ \ x\in\Pi^{\sharp},
\]
with the boundary (\ref{36}) and periodicity (\ref{37}) condition.

In view of (\ref{F9}), (\ref{F3}) the compatibility condition in this problem
becomes%
\begin{equation}
\mp a^{\pm}\zeta^{\prime}+\left\Vert W^{\pm};L^{2}\left(  \varpi\right)
\right\Vert ^{2}=0, \label{P8}%
\end{equation}
where $a^{\pm}=\pm a$ are taken from (\ref{F3}). This formula furnishes the
asymptotic ans\"{a}tze (\ref{P7}) while estimates of remainder are given by
\cite[Ch. 9]{VaTr}.

According to (\ref{P1}) and (\ref{P8}) we have%
\[
\frac{dM_{1}}{d\zeta}\left(  \pm\zeta\right)  =\pm\frac{1}{\zeta^{\prime}}=\pm
a\left\Vert W^{\pm};L^{2}\left(  \varpi\right)  \right\Vert ^{-2}%
\]
and observe that the outgoing Floquet wave in (\ref{44bis}) corresponds to the
point $\left(  +\zeta,M\left(  \zeta\right)  \right)  $ at an ascending arc of
the graph of the \ function $M_{1}$ while the incoming wave to the point
$\left(  -\zeta,M\left(  \zeta\right)  \right)  $ on a descending arc. This
observation in elasticity and acoustics is well-known, cf. \cite[Ch.1]{VoBa},
\cite{na543}, \cite{na569} and \cite{Fliss}, and has two important inferences.
First, the Sommerfeld principle which indicates the direction of propagation
on waves (\ref{44bis}) by their wavenumbers $\pm\zeta$, may become wrong, see,
e.g., the right descending arc in fig. \ref{fig.4}, c. Second, the limiting
absorbtion principle provides the same classification of waves as the
Mandelstam energy principle but may fall through at points of extrema and
inflexion. The latter is the real reason why we have chosen the universal
energy principle. Notice that restriction $M<M^{\sharp}$ has been introduced
in order to unify our notation and to deal with only a couple of Floquet waves.

\section{Detaching asymptotics\label{sect4}}

\subsection{Weighted spaces\label{sect4.1}}

Following \cite{na17} and \cite[\S 3.4]{NaPl}, we study problems
(\ref{A6})-(\ref{A8}) and (\ref{35})-(\ref{37}) in the Kondratiev spaces
$\mathbf{W}_{\beta}^{l}(\Pi^{\sharp})$ obtained by the completion of
$C_{c}^{\infty}(\overline{\Pi}^{\sharp})$ (infinitely differentiable functions
with compact supports) in the weighted norm%
\begin{equation}
||V;\mathbf{W}_{\beta}^{l}(\Pi^{\sharp})||=\left(
{\displaystyle\sum\limits_{k=0}^{l}}
||e^{\beta\left\vert x_{2}\right\vert }\nabla^{k}V;L^{2}(\Pi^{\sharp}%
)||^{2}\right)  ^{1/2}, \label{55}%
\end{equation}
where $l\in\left\{  0,1,2,...\right\}  $ and $\beta\in\mathbb{R}$ are the
smoothness and weight indexes while $\nabla^{k}V$ stands for a family of all
order $k$ derivatives of $V$. This space consists of all functions in
$H_{loc}^{l}\left(  \Pi^{\sharp}\right)  $ with the finite norm (\ref{55}) and
coincides with $H^{l}(\Pi^{\sharp})$ in the case $\beta=0$. However, for
$\beta>0$ $\left(  \beta<0\right)  $ functions in $\mathbf{W}_{\beta}^{l}%
(\Pi^{\sharp})$ decay exponentially as $x_{2}\rightarrow\pm\infty$ (some
growth at infinity is permitted) while the decay/growth rate is governed by
$\beta$. The subspaces $\mathbf{W}_{\beta,per}^{l}(\Pi^{\sharp})$ and
$\mathbf{W}_{\beta,0}^{1}(\Pi^{\sharp};\Gamma^{\sharp})$ are composed from
functions satisfying (\ref{36}) and (\ref{37}), respectively.

Problem (\ref{A1}), (\ref{A2}) in the domain $\Omega^{\sharp}$ which is
infinite in two direction, requires for the weighted space $W_{\beta,\gamma
}^{l\pm}(\Omega^{\sharp})$ obtained by the completion of $C_{c}^{\infty
}(\overline{\Omega^{\sharp}})$ in the norm
\begin{equation}
||v;W_{\beta,\gamma}^{l\pm}(\Omega^{\sharp})||=\left(
{\displaystyle\sum\limits_{k=0}^{l}}
||e^{\beta\left\vert x_{2}\right\vert +\gamma x_{1}^{\pm}}\nabla^{k}%
v;L^{2}(\Omega^{\sharp})||^{2}\right)  ^{1/2} \label{56}%
\end{equation}
depending on two weight indexes and using the variables $x_{1}^{+}=\left\vert
x_{1}\right\vert $ and $x_{1}^{-}=x_{1}$. The subspace $W_{\beta,\gamma
,0}^{l\pm}(\overline{\Omega^{\sharp}})$ takes into account the Dirichlet
condition (\ref{A2}).

To derive the key estimates, we also will use the space $W_{\beta,\gamma
}^{l\pm}\left(  \Omega^{\bullet}\right)  $ in the domain (\ref{4}) in fig.
\ref{fig.1}, b. We underline a crucial difference between norms (\ref{56})
caused by the superscripts $\pm.$ If both $\beta$ and $\gamma$ are positive,
the weight with plus in (\ref{56}) grows exponentially in all directions but
the weight with minus gets the exponential decay when $\left\vert
x_{2}\right\vert <const$ and $x_{1}\rightarrow+\infty$ but still grows in
other radial directions. These properties will allow us to describe
asymptotics of solutions near the open waveguide $\Xi^{+}.$

\subsection{The problem in the perforated strip\label{sect4.2}}

The inhomogeneous problem (\ref{35})-(\ref{37}) is associated with the mapping%
\begin{equation}
\mathbf{W}_{\beta,per}^{2}(\Pi^{\sharp})\cap\mathbf{W}_{\beta,0}^{1}%
(\Pi^{\sharp};\Gamma^{\sharp})=:\mathbf{W}_{\beta}^{\sharp}\ni W\mapsto
\mathbf{A}_{\beta}\left(  \zeta,M\right)  W=-(\left(  \partial_{1}%
+i\zeta\right)  ^{2}+\partial_{2}^{2}+MW)\in\mathbf{W}_{\beta}^{0}(\Pi
^{\sharp}) \label{57}%
\end{equation}
which evidently is continuous for any $\beta\in\mathbb{R}$ but, according to
\cite{na17} and \cite[ Thm.3.4.6, 5.1.4]{NaPl}, is Fredholm if and only if the
segments%
\begin{equation}
\Upsilon_{\beta}^{\sharp}=\left\{  \xi\in\mathbb{C}:\operatorname{Re}\xi
\in\left[  -\frac{\pi}{l_{2}},\frac{\pi}{l_{2}}\right]  ,\ \operatorname{Im}%
\xi=\beta\right\}  \label{58}%
\end{equation}
in the complex plane is free of the $\xi$-spectrum of the quadratic pencil
\cite[Ch. 1]{GoKr}%
\begin{equation}
\mathbb{C}\ni\xi\mapsto(\mathfrak{A}\left(  \xi;\zeta,M\right)  =-\left(
\partial_{1}+i\zeta\right)  ^{2}-\left(  \partial_{2}+i\xi\right)
^{2}-M:H_{per}^{2}\left(  \varpi\right)  \cap H_{0}^{1}\left(  \varpi
;\partial\omega\right)  \rightarrow L^{2}\left(  \varpi\right)  ) \label{401}%
\end{equation}
where $\zeta\in\mathbb{C}$, $M\in\mathbb{R}$ are fixed and $H_{per}^{2}\left(
\varpi\right)  $ is the Sobolev space of functions which are $2l_{j}$-periodic
in $x_{j},$ $j=1,2$ (recall that $l_{1}=1/2$ and compare (\ref{401}) with
(\ref{39})).

\begin{remark}
\label{RemarkSTRIP}Results in \cite{na17} and \cite[Ch. 3 and 5]{NaPl} are
obtained for general boundary-value problems for elliptic systems in smooth
$n$-dimensional domains with periodic outlets to infinity. The presence of the
periodicity conditions (\ref{37}) does not hamper the applicability of those
results since the lateral sides $\left\{  \pm1/2\right\}  \times\mathbb{R}$ of
the strip $\Pi^{\sharp}$ can be identified so that the problem can be posed on
a perforated cylindrical surface in $\mathbb{R}^{3}$. In this way, a literal
repetition of arguments in \cite{na17} proves all assertions in use below.
$\boxtimes.$
\end{remark}

In the case $M<\lambda_{\dag}^{\circ}$ the segment $\Upsilon_{0}^{\sharp}$ is
free of the spectrum of the pencil (\ref{401}) due to definition of the cutoff
value $\lambda_{\dag}^{\circ}$ and formulas (\ref{25}), (\ref{26}). Notice
that the Fredholm property of $\mathbf{A}_{0}\left(  \zeta;M\right)  $ with
any $\zeta\in\left[  -\pi,\pi\right]  $ implies the formula $M\notin\sigma
_{e}^{0}$ and the inclusion $M\in\sigma_{d}^{\sharp}$ means that the subspace
$\ker\mathbf{A}_{0}\left(  \zeta;M\right)  $ is not trivial,
\[
\ker\mathbf{A}_{\beta}\left(  \zeta;M\right)  =\{W\in\mathbf{W}_{\beta
,per}^{2}(\Pi^{\sharp})\cap\mathbf{W}_{\beta,0}^{1}(\Pi^{\sharp}):W\text{
satisfy (\ref{35})}\}.
\]
By inequality (\ref{21}), problem (\ref{11})-(\ref{13}) with $\Lambda\left(
\eta\right)  =M\in\left(  0,\Lambda^{\star}\right)  $ and $\eta\in\mathbb{Y}$
has only trivial solution. Hence, the analytic Fredholm alternative, see,
e.g., \cite[ Thm 1.5.1]{GoKr} shows that, for $M\in(0,M^{\sharp})$, the $\xi
$-spectrum of the pencil $\mathfrak{A}\left(  \cdot;\zeta,M\right)  $,
(\ref{401}), is a countable set of normal eigenvalues without finite
accumulation points. This spectrum is invariant with respect to shifts $\pm
\pi/l_{2}$ along the real axis because the eigenpairs $\left\{  \xi
,U^{0}\left(  x\right)  \right\}  $ and $\left\{  \xi\pm2\pi/l_{2},e^{\mp i\pi
x_{2}/l_{2}}U^{0}\left(  x\right)  \right\}  $ occur simultaneously. Thus,
there exists a positive $\beta^{\sharp}\left(  M\right)  $ such that, for any
$\zeta\in\left[  -\pi,\pi\right]  ,$ the rectangle%
\begin{equation}
\Xi_{\beta^{\sharp}\left(  M\right)  }^{\sharp}=\{\xi\in\mathbb{C}:\left\vert
\operatorname{Re}\xi\right\vert \leq\pi/2l_{2}\text{, }\left\vert
\operatorname{Im}\xi\right\vert <\beta^{\sharp}\left(  M\right)
\}\supset\Upsilon_{0}^{\sharp} \label{59}%
\end{equation}
is free of the $\xi$-spectrum, too. Besides, the theorem on asymptotics, see
\cite[Thm4]{na17}, \cite[Thm 3.4.7 and 5.1.4]{NaPl}, ensures that the kernel
of the operator $\mathbf{A}_{\beta}\left(  \zeta;M\right)  $ is independent of
$\beta\in(-\beta^{\sharp}\left(  M\right)  ,\beta^{\sharp}\left(  M\right)
)$,%
\begin{equation}
\ker\mathbf{A}_{\beta}\left(  \zeta;M\right)  =\ker\mathbf{A}_{0}\left(
\zeta;M\right)  \text{ \ }\forall\zeta\in\left[  -\pi,\pi\right]  . \label{60}%
\end{equation}

\begin{remark}
\label{RemarkDECAY}In view of (\ref{60}) a trapped mode $W^{0}\in
\ker\mathbf{A}_{0}\left(  \zeta;M\right)  \subset H_{per}^{2}(\Pi^{\sharp
})\cap H_{0}^{1}(\Pi^{\sharp};\Gamma^{\sharp})$ falls into $\mathbf{W}%
_{\beta,per}^{2}(\Pi^{\sharp})\cap\mathbf{W}_{\beta,0}^{1}(\Pi^{\sharp}%
;\Gamma^{\sharp})$ and therefore decays exponentially at infinity. Then the
right-hand side $F^{1}$of equation (\ref{41}) belongs to $\mathbf{W}_{\beta
}^{0}(\Pi^{\sharp})$. The formally self-adjoint problem (\ref{41}),
(\ref{36}), (\ref{37}) admits a solution $W^{1}\in\mathbf{W}_{\beta,per}%
^{2}(\Pi^{\sharp})\cap\mathbf{W}_{\beta,0}^{1}(\Pi^{\sharp};\Gamma^{\sharp})$
if and only if
\begin{equation}
\int_{\Pi^{\sharp}}F^{1}\left(  x\right)  \overline{W\left(  x\right)
}dx=0\text{ \ }\forall W\in\ker\mathbf{A}_{-\beta}\left(  0;M\right)  .
\label{61}%
\end{equation}
At the same time, $\ker\mathbf{A}_{-\beta}\left(  0;M\right)  $ is spanned
over the eigenfunction $W^{0}$ because of (\ref{60}) and, thus, (\ref{61})
converts into (\ref{400}). $\boxtimes$
\end{remark}

\subsection{The problem in the periodic perforated plane.\label{sect4.3}}

The composition of the Gelfand transform (\ref{A5}) and the change (\ref{34})
takes the form%
\begin{equation}
v(x)\mapsto W(x;\zeta)=\left(  \mathcal{G}v\right)  \left(  x;\zeta\right)
=\frac{1}{\sqrt{2\pi}}%
{\displaystyle\sum\limits_{\alpha_{1}\in\mathbb{Z}}}
e^{-i\zeta\left(  x_{1}+\alpha_{1}\right)  }v\left(  x_{1}+\alpha_{1}%
,x_{2}\right)  . \label{62}%
\end{equation}
Note that $x\in\Omega^{\sharp}$ on the left but $x\in\Pi^{\sharp}$ on the
right in (\ref{62}). By a direct calculation, cf. \cite[\S 3.4]{NaPl},
transform (\ref{62}), establishes the isometric isomophism%
\begin{equation}
W_{\beta,0}^{0-}(\Omega^{\sharp})\cong L^{2}(\Upsilon_{0};\mathbf{W}_{\beta
}^{0}(\Pi^{\sharp})), \label{63}%
\end{equation}
where%
\begin{equation}
\Upsilon_{\gamma}=\left\{  \zeta\in\mathbb{C}:\left\vert \operatorname{Re}%
\zeta\right\vert \leq\pi,\ \operatorname{Im}\zeta=\gamma\right\}  \label{402}%
\end{equation}
and $L^{2}\left(  \Upsilon_{0};\mathfrak{B}\right)  $ stands for the Lebesgue
space of abstract functions with values in a Banach space $\mathfrak{B}$ and
the norm
\[
\left\Vert W;L^{2}\left(  \Upsilon_{0};\mathfrak{B}\right)  \right\Vert
=\left(  \int_{\Upsilon_{0}}\left\Vert W\left(  \zeta\right)  ;\mathfrak{B}%
\right\Vert ^{2}ds\right)  ^{1/2}.
\]
The change
\[
v\left(  x\right)  \mapsto\mathbf{w}\left(  x\right)  =e^{-\gamma x_{1}%
}v\left(  x\right)
\]
which provides the equivalency of the norms $||\mathbf{w};W_{\beta,0}%
^{0-}(\Omega^{\sharp})||$ and $||v;W_{\beta,\gamma}^{0-}(\Omega^{\sharp})||,$
cf. definition (\ref{56}), passes property (\ref{63}) to the Gelfand transform
(\ref{62}) with $\zeta\in\Upsilon_{\gamma},$ that is, with a complex-valued
dual variable. As a result, we come across the isomorphism, not necessarily
isometric,
\[
W_{\beta,\gamma}^{0}(\Omega^{\sharp})\approx L^{2}(\Upsilon_{\gamma
};\mathbf{W}_{\beta}^{0}(\Pi^{\sharp})).
\]

Let $\lambda\in\left(  0,M_{1}\left(  0\right)  \right)  .$ According to Lemma
\ref{lemmaBBB}, problem (\ref{35})-(\ref{37}) has no trivial solution in
$H^{2}(\Pi^{\sharp})$. Owing to the above-mentioned properties of the $\xi
$-spectrum of $\mathfrak{A}\left(  \cdot,\zeta,\lambda\right)  $ the unique
solvability of the differential equation%
\begin{equation}
-\left(  \partial_{1}+i\zeta\right)  ^{2}W\left(  x,\zeta\right)
-\partial_{2}^{2}W\left(  x,\zeta\right)  -\lambda W\left(  x,\zeta\right)
=F\left(  x,\zeta\right)  ,\text{ \ }x\in\Pi^{\sharp}, \label{64}%
\end{equation}
with the usual conditions (\ref{36}), (\ref{37}) is also kept in
$\mathbf{W}_{\beta,per}^{2}(\Pi^{\sharp})\cap\mathbf{W}_{\beta,0}^{1}%
(\Pi^{\sharp};\Gamma^{\sharp})$ for $F\in\mathbf{W}_{\beta}^{0}(\Pi^{\sharp})$
and $\beta\in(-\beta^{\sharp}\left(  \lambda\right)  ,\beta^{\sharp}\left(
\lambda\right)  )$ where $\beta^{\sharp}\left(  \lambda\right)  >0$ depends on
$\lambda$ and vanishes when $\lambda\rightarrow M_{1}\left(  0\right)  -0.$

Taking $f\in W_{\beta,0}^{0-}(\Omega^{\sharp})$ and applying the Gelfand
transform (\ref{62}), we solve problem (\ref{64}), (\ref{36}), (\ref{37}) with
the right-hand side $F=\mathcal{G}f$ and obtain a unique solution $W\left(
\cdot,\zeta\right)  =-\mathbf{A}_{\beta}\left(  \zeta;M\right)  ^{-1}F\left(
\cdot,\zeta\right)  $ together with the estimate
\[
\left\Vert W\left(  \cdot,\zeta\right)  ;L^{2}(\Upsilon_{0};\mathbf{W}_{\beta
}^{2}(\Pi^{\sharp}))\right\Vert ^{2}\leq c\left\Vert F\left(  \cdot
,\zeta\right)  ;L^{2}(\Upsilon_{0};\mathbf{W}_{\beta}^{0}(\Pi^{\sharp
}))\right\Vert ^{2}\leq C||f;W_{\beta,0}^{0-}(\Omega^{\sharp})||^{2}.
\]
The inverse Gelfand transform acts as follows:%
\begin{equation}
W\left(  x,\zeta\right)  \mapsto v\left(  x\right)  =\frac{1}{\sqrt{2\pi}}%
\int_{\Upsilon_{\gamma}}e^{i\zeta x_{1}}W\left(  x_{1}-\left[  x_{1}\right]
,x_{2},\zeta\right)  d\zeta, \label{inv}%
\end{equation}
see, e.g., \cite{na17} and \cite[\S 3.4]{NaPl}; here, $\left[  t\right]
=\max\left\{  m\in\mathbb{Z}:m\leq t\right\}  $ while $x\in\Omega^{\sharp}$ on
the right in (\ref{inv}). It gives us a solution $v=\mathcal{G}^{-1}W\in
W_{\beta,0}^{2-}(\Omega^{\sharp})\cap W_{\beta,0}^{1-}(\Omega^{\sharp})$ of
problem (\ref{A1}), (\ref{A2}) which meets the estimate%
\[
\left\Vert v;W_{\beta,0}^{2-}(\Omega^{\sharp})\right\Vert ^{2}\leq c\int
_{-\pi}^{\pi}\left\Vert F\left(  \cdot,\zeta\right)  ;\mathbf{W}_{\beta}%
^{0}(\Pi^{\sharp})\right\Vert ^{2}d\zeta\leq C||f;W_{\beta,0}^{0-}%
(\Omega^{\sharp})||^{2}%
\]
and is unique because of the above-mentioned uniqueness of $W\left(
\cdot;\zeta\right)  .$

This standard scheme to solve boundary-value problems in periodic domains
breaks in the case%
\begin{equation}
\lambda\in\lbrack M_{1}\left(  0\right)  ,M^{\sharp}) \label{65}%
\end{equation}
since problem (\ref{35})-(\ref{37}) gets a trapped mode for some $\zeta
\in\left(  -\pi,\pi\right)  .$

\subsection{Asymptotics at infinity\label{sect4.4}}

As was deduced in Sections \ref{sect2.3} and \ref{sect3.1}, problem
(\ref{23}), (\ref{A7}), (\ref{A8}) gains a trapped mode which gives rise to
the Floquet waves (\ref{66}) and (\ref{44bis}) in the homogeneous problem
(\ref{5}), (\ref{6}).

Let $M=\lambda$ in (\ref{65}) be fixed and let $\beta^{\sharp}\left(
M\right)  >0$ be chosen such that rectangle (\ref{59}) in the complex plane
includes just two real points $\zeta=\pm\zeta\left(  M\right)  \in\Upsilon
_{0}^{\sharp}$, recall the notation in Section \ref{sect3.1} and an argument
in Section \ref{sect4.2}. The problem (\ref{35})-(\ref{37}) with $\zeta
=\pm\zeta\left(  M\right)  $ gets the eigenfunction $W^{\pm}\in H_{per}%
^{2}(\Pi^{\sharp})\cap H_{0}^{1}(\Pi^{\sharp};\Gamma^{\sharp})$, cf.
(\ref{44}), which, according to Remark \ref{RemarkDECAY}, falls into
$\mathbf{W}_{\beta,per}^{2}(\Pi^{\sharp})\cap\mathbf{W}_{\beta,0}^{1}%
(\Pi^{\sharp};\Gamma^{\sharp})$ with some $\beta>0$ and, therefore, $\zeta
=\pm\zeta\left(  M\right)  $ are real eigenvalues of the quadratic pencil%
\begin{equation}
\mathbb{C}\ni\zeta\mapsto(\mathfrak{A}_{\beta}^{\sharp}\left(  \zeta;M\right)
:\mathbf{W}_{\beta,per}^{2}(\Pi^{\sharp})\cap\mathbf{W}_{\beta,0}^{1}%
(\Pi^{\sharp};\Gamma^{\sharp})\rightarrow\mathbf{W}_{\beta}^{0}(\Pi^{\sharp
})). \label{67}%
\end{equation}
Moreover, there is no other eigenvalue of (\ref{67}) in the segment
$\Upsilon_{0},$ (\ref{402}), and we can fix $\gamma>0$ such that
$\Upsilon_{\pm\gamma}$ is free of the $\zeta$-spectrum of the pencil. As a
result, repeating an argumentation in the end of Section \ref{sect4.3} with
the replacement $\Upsilon_{0}\mapsto\Upsilon_{\pm}$ and using the Gelfand
transform with complex dual variable deliver two solutions%
\begin{equation}
v^{\pm}\left(  x\right)  =\frac{1}{\sqrt{2\pi}}\int_{\pm\gamma}e^{i\zeta
x_{1}}\mathbf{A}_{\beta}\left(  \zeta;M\right)  ^{-1}F\left(  x_{1}-\left[
x_{1}\right]  ,x_{2};\zeta\right)  d\zeta\label{68}%
\end{equation}
of problem (\ref{A1}), (\ref{A2}) with the right-hand side%
\begin{equation}
f\in W_{\beta,\gamma}^{0-}(\Omega^{\sharp})\cap W_{\beta,-\gamma}^{0-}%
(\Omega^{\sharp})=W_{\beta,\gamma}^{0+}(\Omega^{\sharp}). \label{69}%
\end{equation}
Due to the exponential decay of the function $f$ as $x_{1}\rightarrow\pm
\infty$, compare (\ref{69}) and (\ref{56}), the Gelfand transform
$F=\mathcal{G}f$ is an analytic abstract function in the variable $\zeta\in
\Xi_{\gamma}$ with values in $\mathbf{W}_{\beta}^{0}(\Pi^{\sharp})$; moreover,
$F$ is $2\pi$-periodic in $\operatorname{Re}\zeta$ and continuous up to the
boundary of the open rectangle%
\[
\Xi_{\gamma}=\left\{  \zeta:\left\vert \operatorname{Re}\zeta\right\vert
<\pi,\ \left\vert \operatorname{Im}\zeta\right\vert <\gamma\right\}  .
\]
Notice that each of solutions (\ref{68}) is unique in its own class
$W_{\beta,\pm\gamma}^{2-}(\Omega^{\sharp})\cap W_{\beta,\pm\gamma,0}%
^{1-}(\Omega^{\sharp}).$

We are in position to apply the Cauchy residue theorem to the contour integral%
\begin{equation}
\frac{1}{\sqrt{2\pi}}%
{\displaystyle\oint_{\partial\Xi_{\gamma}}}
e^{izx_{1}}\mathbf{A}_{\beta}\left(  \zeta;M\right)  ^{-1}F\left(  x-\left[
x_{1}\right]  ,x_{2};\zeta\right)  d\zeta. \label{70}%
\end{equation}
We observe that in view of $2\pi$-periodicity along the real axis, the
integrals along the vertical sides of the rectangle cancel each other while
the sum of the integrals along the horizontal sides equals the difference
$v^{-}\left(  x\right)  -$ $v^{+}\left(  x\right)  $ (in both cases direction
of integration was taken into account). At the same time, the contour integral
turns in the sum of residuals which are to be computed according to the
formula%
\begin{equation}
\mathbf{A}_{\beta}\left(  \zeta;M\right)  ^{-1}=\left(  \zeta-\zeta_{\pm}%
^{\pm}\right)  ^{-1}W^{\pm}R^{\pm}+\mathbf{R}^{\pm}\left(  \zeta;M\right)
,\text{ \ }\zeta\in\mathbb{B}_{\rho}\left(  \zeta^{\pm}\right)  , \label{71}%
\end{equation}
where $\mathbb{B}_{\rho}\left(  \zeta^{\pm}\right)  =\left\{  \zeta
\in\mathbb{C}:\left\vert \zeta-\zeta^{\pm}\right\vert <\rho\right\}  $ is a
disk of a small radius $\rho>0,$ $\zeta^{\pm}=\pm\zeta$ and $W^{\pm}$
respectively are simple eigenvalues and the corresponding eigenvectors of
pencil (\ref{39}), see Section \ref{sect3.1}, $R^{\pm}$ is a continuous
functional in $\mathbf{W}_{\beta}^{0}(\Pi^{\sharp})$ and $\mathbf{R}\left(
\cdot;M\right)  $ is analytic in $\mathbb{B}_{\rho}\left(  \zeta^{\pm}\right)
$. As a result, we conclude the representation
\begin{equation}
v^{-}\left(  x\right)  =v^{+}\left(  x\right)  +a_{+}w^{+}\left(  x\right)
+a_{-}w^{-}\left(  x\right)  \label{72}%
\end{equation}
where $w^{\pm}$ are the Floquet waves (\ref{44bis}) and $a_{\pm}$ are some
coefficients satisfying the estimate%
\begin{equation}
\left\vert a_{+}\right\vert +\left\vert a_{-}\right\vert \leq c_{\gamma
}||f;W_{\beta,\gamma}^{0+}(\Omega^{\sharp})||. \label{73}%
\end{equation}

Recalling the notation (\ref{42}), (\ref{42N}) and (\ref{66}), we derive the
same formulas (\ref{72}) and (\ref{73}) also in the case $\lambda=M_{1}\left(
0\right)  $ when the eigenvalue $\zeta=0$ of the pencil $\mathfrak{A}\left(
\cdot;\lambda\right)  $ is of algebraic multiplicity $2$ and generates the
Jordan chain $\left\{  W,W^{1}\right\}  $. We only mention that the new
resolvent (\ref{71}) gains a pole of degree $2$ at the point $\zeta=0$.

\begin{theorem}
\label{TheoremASY}For $\lambda\in\lbrack M_{1}\left(  0\right)  ,M^{\sharp})$,
problem (\ref{A1}), (\ref{A2}) in $\Omega^{\sharp},$ see (\ref{A3}) and fig.
\ref{fig.3}, a, with the right-hand side (\ref{69}) has two solutions $v^{\pm
}\in W_{\beta,\pm\gamma}^{2-}(\Omega^{\sharp})\cup W_{\beta,\pm\gamma,0}%
^{1-}(\Omega^{\sharp})$ which are given in (\ref{68}) and are related by the
asymptotic formula (\ref{72}), where $w^{\pm}$ are the Floquet waves
(\ref{66}) or (\ref{44bis}) and the coefficients $a_{\pm}$ enjoy estimate
(\ref{73}).
\end{theorem}

In the next section we will interpret (\ref{72}) as an asymptotic
decomposition of the growing solution $v^{-}\left(  x\right)  $ with the
decaying remainder $v^{+}\left(  x\right)  $ when $x_{1}\rightarrow+\infty.$

\section{Solvability of the problem with radiation conditions\label{sect5}}

\subsection{The localization estimates\label{sect5.1}}

The integral identity
\begin{equation}
\left(  \nabla u,\nabla v\right)  _{\Omega^{\bullet}}-\lambda\left(
u,v\right)  _{\Omega^{\bullet}}=f\left(  v\right)  \text{ \ }\forall v\in
W_{\beta,\gamma,0}^{1+}\left(  \Omega^{\bullet}\right)  \label{74}%
\end{equation}
serves for the inhomogeneous problem (\ref{5}), (\ref{6}) in the weighted
space $W_{-\beta,-\gamma,0}^{1+}\left(  \Omega^{\bullet}\right)  \ni u$.
According to definition (\ref{56}) all terms in (\ref{74}) are defined
properly if $f\in W_{\beta,\gamma,0}^{1+}\left(  \Omega^{\bullet}\right)
^{\ast}$ is an (anti)linear functional in $W_{\beta,\gamma,0}^{1+}\left(
\Omega^{\bullet}\right)  $ and $\left(  \ ,\ \right)  _{\Omega^{\bullet}}$ is
understood as an extension of the scalar product in $L^{2}\left(
\Omega^{\bullet}\right)  $ up to the duality between $\mathcal{L}%
_{-\beta,-\gamma}\left(  \Omega^{\bullet}\right)  $ and $\mathcal{L}%
_{\beta,\gamma}\left(  \Omega^{\bullet}\right)  $. Here, $\mathcal{L}%
_{\beta,\gamma}\left(  \Omega^{\bullet}\right)  $ is a weighted Lebesgue space
with the norm%
\begin{equation}
\left\Vert f;\mathcal{L}_{\beta,\gamma}\left(  \Omega^{\bullet}\right)
\right\Vert =||e^{\beta\left\vert x_{1}\right\vert +\gamma\left\vert
x_{2}\right\vert }f;L^{2}\left(  \Omega^{\bullet}\right)  ||. \label{75}%
\end{equation}
The weak formulation (\ref{74}) of the problem in $\Omega^{\bullet}$ generates
the continuous mapping
\begin{equation}
W_{-\beta,-\gamma,0}^{1+}\left(  \Omega^{\bullet}\right)  \ni u\mapsto
\mathcal{A}_{_{-\beta,-\gamma}}\left(  M\right)  u=f\in W_{\beta,\gamma
,0}^{1+}\left(  \Omega^{\bullet}\right)  ^{\ast} \label{76}%
\end{equation}
while $\mathcal{A}_{_{\beta,\gamma}}\left(  M\right)  $ is adjoint for
$\mathcal{A}_{_{-\beta,-\gamma}}\left(  M\right)  .$

The following assertion provides the key localization estimate which
demonstrates that a growing solution of the problem with a decaying right-hand
side gets the decay property outside a sectorial neighborhood of the inclusion
(\ref{10}).

\begin{lemma}
\label{lemmaM1}Let%
\begin{equation}
\lambda<\Lambda^{\star},\text{ \ }\beta,\gamma>0\text{ \ and \ \ }%
\lambda+\left(  \beta+\gamma\right)  ^{2}<\Lambda^{\star}. \label{77}%
\end{equation}
Then a solution $u\in W_{-\beta,-\gamma,0}^{1+}\left(  \Omega^{\bullet
}\right)  $ of problem (\ref{74}) with the right-hand side%
\begin{equation}
f\left(  v\right)  =\left(  f,v\right)  _{\Omega^{\bullet}},\text{ \ }%
f\in\mathcal{L}_{\beta,\gamma}^{2}\left(  \Omega^{\bullet}\right)  \label{777}%
\end{equation}
belongs to the space $W_{\beta,-\gamma,0}^{1-}\left(  \Omega^{\bullet}\right)
$ and obeys the estimate%
\begin{equation}
||u;W_{\beta,-\gamma,0}^{1-}\left(  \Omega^{\bullet}\right)  ||\leq
c||f;\mathcal{L}_{\beta,\gamma}^{2}\left(  \Omega^{\bullet}\right)
||+||u;W_{-\beta,-\gamma,0}^{1+}\left(  \Omega^{\bullet}\right)  || \label{78}%
\end{equation}
where the factor $c$ depends on $\lambda$ and $\beta,$ $\gamma$ but is
independent of $f$ and $u$.
\end{lemma}

\textbf{Proof.} Borrowing a trick from \cite{CaDuNaSIAM}, we introduce the
continuous function $\mathcal{R}_{R}\left(  x\right)  =\mathcal{R}_{R1}\left(
x_{1}\right)  \times\mathcal{R}_{R2}\left(  x_{2}\right)  $ where $R>0$ is a
big parameter and%
\begin{equation}
\mathcal{R}_{R1}\left(  x\right)  =\left\{
\begin{array}
[c]{c}%
e^{-\gamma x_{1}},\ x_{1}\geq-R,\\
e^{2\gamma R}e^{-\gamma\left\vert x_{1}\right\vert },\ x_{1}\leq-R,
\end{array}
\right.  \ \ \ \mathcal{R}_{R2}\left(  x\right)  =\left\{
\begin{array}
[c]{c}%
e^{\beta\left\vert x_{2}\right\vert },\ \left\vert x_{2}\right\vert \leq R,\\
e^{2\beta R}e^{-\beta\left\vert x_{2}\right\vert },\ \left\vert x_{2}%
\right\vert \geq R.
\end{array}
\right.  \label{79}%
\end{equation}
We set $u_{R}=\mathcal{R}_{R}u,$ $v_{R}=\mathcal{R}_{R}u_{R}$ and observe that
$u_{R}\in H_{0}^{1}\left(  \Omega^{\bullet}\right)  ,$ $v_{R}\in
W_{-\beta,-\gamma,0}^{1+}\left(  \Omega^{\bullet}\right)  $ because%
\begin{equation}
\mathcal{R}_{R}\left(  x\right)  \leq c_{R}e^{-\beta\left\vert x_{2}%
\right\vert -\gamma\left\vert x_{1}\right\vert },\text{ \ \ }\left\vert
\nabla\mathcal{R}_{R}\left(  x\right)  \right\vert \leq\left(  \beta
+\gamma\right)  \mathcal{R}_{R}\left(  x\right)  . \label{80}%
\end{equation}
Inserting $v_{R}$ as a test function into (\ref{74}) and performing simple
algebraic transformations yield%
\begin{align*}
\left(  \mathcal{R}_{R}f,u_{R}\right)  _{\Omega^{\bullet}}  &  =\left\Vert
\nabla u_{R};L^{2}\left(  \Omega^{\bullet}\right)  \right\Vert ^{2}%
-\lambda\left\Vert u_{R};L^{2}\left(  \Omega^{\bullet}\right)  \right\Vert
^{2}\\
&  -\left\Vert u_{R}\mathcal{R}_{R}^{-1}\nabla\mathcal{R}_{R};L^{2}\left(
\Omega^{\bullet}\right)  \right\Vert ^{2}+\left(  \left(  \nabla u_{R}%
,u_{R}\mathcal{R}_{R}^{-1}\nabla\mathcal{R}_{R}\right)  _{\Omega^{\bullet}%
}-\left(  u_{R}\mathcal{R}_{R}^{-1}\nabla\mathcal{R}_{R},\nabla u_{R}\right)
_{\Omega^{\bullet}}\right)  .
\end{align*}
The last difference in brackets is pure imaginary. Hence,%
\begin{align}
&  \left\Vert \nabla u_{R};L^{2}\left(  \Omega^{\bullet}\setminus\Xi
^{+}\right)  \right\Vert ^{2}-\lambda\left\Vert u_{R};L^{2}\left(
\Omega^{\bullet}\setminus\Xi^{+}\right)  \right\Vert ^{2}-\left\Vert
u_{R}\mathcal{R}_{R}^{-1}\nabla\mathcal{R}_{R};L^{2}\left(  \Omega^{\bullet
}\right)  \right\Vert ^{2}\label{81}\\
&  =\operatorname{Re}\left(  \mathcal{R}_{R}f,u_{R}\right)  _{\Omega^{\bullet
}}-\left\Vert \nabla u_{R};L^{2}\left(  \Xi^{+}\right)  \right\Vert
^{2}+\lambda\left\Vert u_{R};L^{2}\left(  \Xi^{+}\right)  \right\Vert
^{2}+\left\Vert u_{R}\mathcal{R}_{R}^{-1}\nabla\mathcal{R}_{R};L^{2}\left(
\Xi^{+}\right)  \right\Vert ^{2}.\nonumber
\end{align}
By (\ref{79}), we have $\mathcal{R}_{R}\left(  x\right)  =e^{-\gamma x_{1}}$
in $\Xi^{+}$ for a big $R$ and $\mathcal{R}_{R}\left(  x\right)  \leq
e^{\gamma\left\vert x_{1}\right\vert +\beta\left\vert x_{2}\right\vert }$ in
$\Omega^{\bullet}$. Hence, the right-hand side of (\ref{81}) does not exceed
the expression%
\begin{equation}
c\left(  \left\Vert f;\mathcal{L}_{\beta,\gamma}^{2}\left(  \Omega^{\bullet
}\right)  \right\Vert \left\Vert u_{R};L^{2}\left(  \Omega^{\bullet}\right)
\right\Vert +||u;W_{-\beta,-\gamma}^{1+}\left(  \Xi^{+}\right)  ||^{2}\right)
. \label{82}%
\end{equation}
According to formula (\ref{4}) the set $\Omega^{\bullet}\setminus\Xi^{+}$
consists of the cells $\varpi\left(  \alpha\right)  $, (\ref{3}), with
$\alpha\in\mathcal{Z}=\mathbb{Z}^{2}\setminus\left\{  \alpha:\alpha_{1}%
\geq1,\ \alpha_{2}\in\left[  1,J\right]  \right\}  $ while the Friedrichs
inequality (\ref{21}) leads to the relation%
\[
\left\Vert \nabla u_{R};L^{2}\left(  \Omega^{\bullet}\setminus\Xi^{+}\right)
\right\Vert ^{2}=%
{\displaystyle\sum\nolimits_{\alpha\in\mathcal{Z}}}
\left\Vert \nabla u_{R};L^{2}\left(  \varpi\left(  \alpha\right)  \right)
\right\Vert ^{2}\geq\Lambda^{\star}%
{\displaystyle\sum\nolimits_{\alpha\in\mathcal{Z}}}
\left\Vert u_{R};L^{2}\left(  \varpi\left(  \alpha\right)  \right)
\right\Vert ^{2}=\Lambda^{\star}\left\Vert u_{R};L^{2}\left(  \Omega^{\bullet
}\setminus\Xi^{+}\right)  \right\Vert ^{2}.
\]
Taking the last formulas in (\ref{77}) and (\ref{80}) into account, we find
some $\delta>0$ such that the left-hand side of (\ref{81}) is bigger than%
\begin{align}
&  \delta\left\Vert \nabla u_{R};L^{2}\left(  \Omega^{\bullet}\setminus\Xi
^{+}\right)  \right\Vert ^{2}+\delta\left\Vert u_{R};L^{2}\left(
\Omega^{\bullet}\setminus\Xi^{+}\right)  \right\Vert ^{2}\label{83}\\
&  =\delta\left\Vert \mathcal{R}_{R}\nabla u+u\nabla\mathcal{R}_{R}%
;L^{2}\left(  \Omega^{\bullet}\setminus\Xi^{+}\right)  \right\Vert ^{2}%
+\delta\left\Vert \mathcal{R}_{R}u;L^{2}\left(  \Omega^{\bullet}\setminus
\Xi^{+}\right)  \right\Vert ^{2}\nonumber\\
&  \geq\delta\tau\left\Vert \mathcal{R}_{R}\nabla u;L^{2}\left(
\Omega^{\bullet}\setminus\Xi^{+}\right)  \right\Vert ^{2}+\delta\left(
1-\tau\frac{\left(  \beta+\gamma\right)  ^{2}}{1-\tau}\right)  \left\Vert
\mathcal{R}_{R}u;L^{2}\left(  \Omega^{\bullet}\setminus\Xi^{+}\right)
\right\Vert ^{2}.\nonumber
\end{align}
Here, we applied the second relation in (\ref{80}) and the simple formula
$\left(  \alpha+\beta\right)  ^{2}\geq\tau a^{2}-\tau\left(  1-\tau\right)
^{-1}b^{2}$ with $\tau\in\left(  0,1\right)  $ and $\tau\left(  1-\tau\right)
^{-1}\left(  \beta+\gamma\right)  ^{2}\geq1/2$. We add the expression%
\[
\left\Vert \mathcal{R}_{R}\nabla u;L^{2}\left(  \Xi^{+}\right)  \right\Vert
^{2}+\left\Vert \mathcal{R}_{R}u;L^{2}\left(  \Xi^{+}\right)  \right\Vert
^{2}=\left\Vert e^{-\gamma x_{1}}\nabla u;L^{2}\left(  \Xi^{+}\right)
\right\Vert ^{2}+\left\Vert e^{-\gamma x_{1}}u;L^{2}\left(  \Xi^{+}\right)
\right\Vert ^{2}%
\]
to both sides of equality (\ref{81}) and estimate its fragments by means of
(\ref{82}) and (\ref{83}). As a result, we obtain the inequality%
\begin{equation}
\left\Vert \mathcal{R}_{R}\nabla u;L^{2}\left(  \Omega^{\bullet}\right)
\right\Vert ^{2}+\left\Vert \mathcal{R}_{R}u;L^{2}\left(  \Omega^{\bullet
}\right)  \right\Vert ^{2}\leq C\left(  \left\Vert f;\mathcal{L}_{\beta
,\gamma}^{2}\left(  \Omega^{\bullet}\right)  \right\Vert ^{2}+||u;W_{-\beta
,-\gamma}^{1+}\left(  \Xi^{+}\right)  ||^{2}\right)  \label{84}%
\end{equation}
where $C$ is independent of $f$, $u$ and $R$. Comparing (\ref{79}) and
(\ref{56}), we see that the left-hand side of (\ref{84}) exceeds
$||u;W_{\beta,-\gamma}^{1-}\left(  \Omega_{R}^{\bullet}\right)  ||^{2}$ where
$\Omega_{R}^{\bullet}=\{x\in\Omega^{\bullet}:\left\vert x_{j}\right\vert <R,$
$j=1,2\}.$ Thus, the limit passage $R\rightarrow\infty$ in (\ref{84}) provides
estimate (\ref{78}) and, therefore, the inclusion $u\in W_{\beta,-\gamma
,0}^{1-}\left(  \Omega^{\bullet}\right)  $ is valid, too. $\boxtimes$

\bigskip

To apply in the next section Theorem \ref{TheoremASY} an asymptotics, we prove
the following lemma which lifts the smoothness of the weak solution.

\begin{lemma}
\label{lemmaM2}Under the condition of Lemma \ref{lemmaM1}, the solution $u$ of
problem (\ref{74}) falls into $\mathbf{W}_{\beta,-\gamma}^{2-}\left(
\Omega^{\bullet}\right)  $ and fulfils the estimate%
\begin{equation}
||u;W_{\beta,-\gamma}^{2-}\left(  \Omega^{\bullet}\right)  ||\leq
c||f;\mathcal{L}_{\beta,\gamma}^{2}\left(  \Omega^{\bullet}\right)
||+||u;W_{\beta,-\gamma}^{0-}\left(  \Omega^{\bullet}\right)  ||. \label{85}%
\end{equation}

\end{lemma}

\textbf{Proof.} For $\alpha\in\mathbb{Z}^{2}$ and $p=0,1,$ we determine the
subdomains $\Omega_{p}^{\bullet}\left(  \alpha\right)  =\left\{  x\in
\Omega^{\bullet}:\left\vert x_{j}-2\alpha_{j}l_{j}\right\vert <\left(
1+p/2\right)  l_{j}\right\}  $ and apply local estimates \cite{ADN1} of
solutions to the Dirichlet problem for the inhomogeneous Helmholtz equation,
namely%
\begin{equation}
\left\Vert \nabla^{2}u;L^{2}\left(  \Omega_{0}^{\bullet}\left(  \alpha\right)
\right)  \right\Vert ^{2}\leq c\left(  \left\Vert f;L^{2}\left(  \Omega
_{1}^{\bullet}\left(  \alpha\right)  \right)  \right\Vert ^{2}+\left\Vert
u;L^{2}\left(  \Omega_{1}^{\bullet}\left(  \alpha\right)  \right)  \right\Vert
^{2}\right)  . \label{86}%
\end{equation}
Owing to the periodic structure of $\Omega^{\bullet},$ we detect only finite
number of homothetically different couples $\Omega_{0}^{\bullet}\left(
\alpha\right)  \subset\Omega_{1}^{\bullet}\left(  \alpha\right)  $ and
therefore can fix the factor $c$ in (\ref{86}) independent of $\alpha$ and, of
course, of $f$ and $u.$ Moreover,
\[
0<c_{\beta,\gamma}\leq\left(  \sup\nolimits_{x\in\Omega_{0}^{\bullet}\left(
\alpha\right)  }e^{\beta\left\vert x_{2}\right\vert -\gamma x_{1}}\right)
^{-1}\inf\nolimits_{x\in\Omega_{1}^{\bullet}\left(  \alpha\right)  }%
e^{\beta\left\vert x_{2}\right\vert -\gamma x_{1}}\leq C_{\beta,\gamma}.
\]
Inserting the exponential weights inside norms in (\ref{86}) and summing in
$\alpha\in\mathbb{Z}^{2}$ yield%
\[
||e^{\beta\left\vert x_{2}\right\vert -\gamma x_{1}}\nabla^{2}u;L^{2}\left(
\Omega^{\bullet}\right)  ||^{2}\leq C\left(  ||e^{\beta\left\vert
x_{2}\right\vert -\gamma x_{1}}f;L^{2}\left(  \Omega^{\bullet}\right)
||^{2}+||e^{\beta\left\vert x_{2}\right\vert -\gamma x_{1}}u;L^{2}\left(
\Omega^{\bullet}\right)  ||^{2}\right)  .
\]
Enlarging the weight of $f$ and taking (\ref{78}) into account lead us to
(\ref{85}). $\boxtimes$

\subsection{The problem with radiation conditions\label{sect5.2}.}

Let%
\begin{equation}
u\in W_{-\beta,-\gamma}^{2+}\left(  \Omega^{\bullet}\right)  \cap
W_{-\beta,-\gamma,0}^{1+}\left(  \Omega^{\bullet}\right)  \label{87}%
\end{equation}
be a solution of the problem%
\begin{equation}
-\Delta u\left(  x\right)  -\lambda u\left(  x\right)  =f\left(  x\right)
,\ x\in\Omega^{\bullet},\ \ \ \ u(x)=0,\text{ \ }x\in\partial\Omega^{\bullet},
\label{88}%
\end{equation}
cf. (\ref{5}), (\ref{6}), with the right-hand side%
\begin{equation}
f\in\mathcal{L}_{\beta,\gamma}^{2}\left(  \Omega^{\bullet}\right)  \label{888}%
\end{equation}
while the spectral parameter and the weight indexes satisfy (\ref{77}), i.e.,
$\beta>0$ and $\gamma>0$ are sufficiently small. We multiply the solution
(\ref{87}) which, by Lemmas \ref{lemmaM1} and \ref{lemmaM2}, belongs to
$W_{\beta,-\gamma}^{2-}\left(  \Omega^{\bullet}\right)  $ with the cut-off
function%
\begin{equation}
\chi\in C^{\infty}\left(  \mathbb{R}\right)  ,\text{ \ }\chi\left(
x_{1}\right)  =1\text{ \ for }x_{1}\geq2,\text{ \ }\chi\left(  x_{2}\right)
=0\ \text{\ \ for }x_{1}\leq1 \label{89}%
\end{equation}
and arrive at problem (\ref{A1}), (\ref{A2}) in $\Omega^{\sharp}$ for
$u^{\chi}=\chi u$ with the new right-hand side%
\[
f^{\chi}=\chi f+\left[  \Delta,\chi\right]  u\in\mathcal{L}_{\beta,\gamma}%
^{2}(\Omega^{\sharp}).
\]
The inclusion holds true because the commutator $\left[  \Delta,\chi\right]
u=2\nabla u\cdot\nabla\chi+u\Delta\chi$ has a support in the closed perforated
strip $\overline{\Pi^{\sharp}\left(  1\right)  }=\{x\in\overline
{\Omega^{\bullet}}:1\leq x_{1}\leq2\}\subset\overline{\Omega^{\sharp}}$ where
the multipliers $e^{\gamma x_{1}^{\pm}}$ do not affect weights in norms
(\ref{56}). Hence, we can apply Theorem \ref{TheoremASY} and conclude the
representation (\ref{72}) for $u^{\chi},$ namely%
\begin{equation}
u^{\chi}\left(  x\right)  =a_{+}w^{+}\left(  x\right)  +a_{-}w^{-}\left(
x\right)  +\widetilde{u}^{\chi}\left(  x\right)  \label{90}%
\end{equation}
with the remainder $\widetilde{u}^{\chi}\in W_{\beta,\gamma}^{2-}%
(\Omega^{\sharp})\cap W_{\beta,\gamma,0}^{1-}(\Omega^{\sharp}).$ Observing
that $e^{\gamma x_{1}^{-}}=e^{\gamma x_{1}^{+}}$ for $x_{1}>1$ and $e^{-\gamma
x_{1}^{-}}\geq e^{-2\gamma}e^{\gamma x_{1}^{+}}$ for $x_{1}<1,$ we have%
\[
\chi\widetilde{u}^{\chi}\in W_{\beta,\gamma}^{2+}\left(  \Omega^{\bullet
}\right)  ,\text{ \ }\left(  1-\chi^{2}\right)  u\in W_{\beta,\gamma}%
^{2+}\left(  \Omega^{\bullet}\right)  .
\]
Then we multiply (\ref{90}) by $\chi$, add $\left(  1-\chi^{2}\right)  u$ to
the result and derive the following representation of solution (\ref{87}):%
\begin{equation}
u\left(  x\right)  =\chi\left(  x_{1}\right)  \left(  a^{+}w^{+}\left(
x\right)  +a^{-}w^{-}\left(  x\right)  \right)  +\widetilde{u}\left(
x\right)  \label{91}%
\end{equation}
together with the estimate%
\begin{equation}
\left\vert a^{+}\right\vert +\left\vert a^{-}\right\vert +||\widetilde
{u};W_{\beta,\gamma}^{2+}\left(  \Omega^{\bullet}\right)  ||\leq c\left(
||f;\mathcal{L}_{\beta,\gamma}^{2}\left(  \Omega^{\bullet}\right)
||+||u;W_{-\beta,-\gamma}^{2+}\left(  \Omega^{\bullet}\right)  ||\right)  .
\label{92}%
\end{equation}
To conclude with (\ref{92}), it should be mentioned that all inclusions
written above are accompanied with estimates of the corresponding norms by the
same majorants as in (\ref{92}).

\begin{theorem}
\label{Theorem ASYFin}Let $\lambda$ and $\beta,\gamma$ satisfy (\ref{77}). A
solution (\ref{87}) of problem (\ref{88}) with the right-hand side (\ref{888})
takes the form (\ref{91}) and estimate (\ref{92}) is valid.
\end{theorem}

In the case $a^{-}=0$ we say that solution (\ref{91}) satisfies
the\ Mandelstam radiation conditions. Indeed, it loses the incoming wave
$w^{-}$ and differs from the outgoing wave $a^{+}\chi w^{+}$ localized near
the semi-infinite inclusion $\Xi^{+},$ by a function with the exponential
decay in all directions.

\subsection{Solvability of the problem with the radiation
condition\label{sect5.3}}

We proceed with the following assertion.

\begin{theorem}
\label{TheoremInd}Let $\lambda$ and $\beta,$ $\gamma$ meet the conditions%
\begin{equation}
\lambda\in\lbrack M_{1}\left(  0\right)  ,M^{\sharp})\text{ and }\beta
,\gamma>0,\text{ \ }\lambda+\beta^{2}+\gamma^{2}<\Lambda^{\star}, \label{Cbg}%
\end{equation}
cf. (\ref{77}). The operators of problem (\ref{88})%
\begin{equation}
\mathcal{A}_{\pm\beta,\pm\gamma}^{2}:W_{\pm\beta,\pm\gamma}^{2+}\left(
\Omega^{\bullet}\right)  \cap W_{\pm\beta,\pm\gamma,0}^{1+}\left(
\Omega^{\bullet}\right)  \mapsto\mathcal{L}_{\pm\beta,\pm\gamma}^{2}\left(
\Omega^{\bullet}\right)  \label{93}%
\end{equation}
are Fredholm and their indexes are as follows:%
\begin{equation}
\mathrm{Ind}\mathcal{A}_{\pm\beta,\pm\gamma}^{2}=\dim\ker\mathcal{A}_{\pm
\beta,\pm\gamma}^{2}-\mathrm{coker}\mathcal{A}_{\pm\beta,\pm\gamma}^{2}=\mp1.
\label{94}%
\end{equation}

\end{theorem}

\textbf{Proof. }To verify the Fredholm property, we follow a scheme proposed
in \cite{CaNaTa} and construct a (right) parametrix $\mathcal{R}_{\pm\beta
,\pm\gamma}^{2}$ for operator (\ref{93}), that is, a continuous mapping
$\mathcal{L}_{\pm\beta,\pm\gamma}^{2}\left(  \Omega^{\bullet}\right)  \mapsto
W_{\pm\beta,\pm\gamma}^{2+}\left(  \Omega^{\bullet}\right)  \cap W_{\pm
\beta,\pm\gamma,0}^{1+}\left(  \Omega^{\bullet}\right)  $ such that
$\mathcal{A}_{\pm\beta,\pm\gamma}^{2}\mathcal{R}_{\pm\beta,\pm\gamma}^{2}-Id$
is a compact operator in $\mathcal{L}_{\pm\beta,\pm\gamma}^{2}\left(
\Omega^{\bullet}\right)  .$ Let us outline this scheme with minor
modifications. First of all, thanks to the "lifting procedure" in our proof of
Lemma \ref{lemmaM2}, we may consider the weak formulation of problem
(\ref{88}) in the space $W_{\pm\beta,\pm\gamma,0}^{1+}\left(  \Omega^{\bullet
}\right)  $, namely%
\begin{equation}
\left(  \triangledown u^{\bullet},\triangledown v^{\bullet}\right)
_{\Omega^{\bullet}}-\lambda\left(  u^{\bullet},v^{\bullet}\right)
_{\Omega^{\bullet}}=f^{\bullet}\left(  v^{\bullet}\right)  \text{ \ }\forall
v^{\bullet}\in W_{\mp\beta,\mp\gamma}^{1+}\left(  \Omega^{\bullet}\right)
\label{95}%
\end{equation}
and the corresponding operator%
\begin{equation}
\mathcal{A}_{\pm\beta,\pm\gamma}^{1\bullet}:W_{\pm\beta,\pm\gamma,0}%
^{1+}\left(  \Omega^{\bullet}\right)  \mapsto W_{\mp\beta,\pm\gamma,0}%
^{1+}\left(  \Omega^{\bullet}\right)  ^{\ast}. \label{96}%
\end{equation}
We also will need formula (\ref{95}) and (\ref{96}) with the change
$\bullet\mapsto\circ$ of the superscript. Taking $f^{\bullet}\in
\mathcal{L}_{\pm\beta,\pm\gamma}^{2}\left(  \Omega^{\bullet}\right)  $, we
annull this function on the foreign inclusion by setting $f^{\circ
}=Xf^{\bullet}\in\mathcal{L}_{\pm\beta,\pm\gamma}^{2}\left(  \Omega^{\circ
}\right)  $ where $\mathcal{X}\in C^{\infty}\left(  \mathbb{R}^{2}\right)  $
is a cut-off function such that
\begin{align*}
\mathcal{X}\left(  x\right)   &  =0\text{ \ for }x\in\Xi^{+}\text{, see
(\ref{10}), and}\\
\mathcal{X}\left(  x\right)   &  =1\text{ \ for }x\in\mathbb{R}^{2}%
\setminus\Xi^{\circ},\ \ \Xi^{\circ}=\{x:x_{1}>1,\text{ }0<x_{2}<2\left(
J+1\right)  l_{2}\}.
\end{align*}
Then we perform the substitutions%
\begin{align}
u^{\circ}\left(  x\right)   &  =e^{\mp\beta\left\vert x_{2}\right\vert
\mp\gamma\left\vert x_{1}\right\vert }u^{\bullet}\left(  x\right)  ,\text{
\ \ }v^{\circ}\left(  x\right)  =e^{\pm\beta\left\vert x_{2}\right\vert
\pm\gamma\left\vert x_{1}\right\vert }v^{\bullet}\left(  x\right)
,\label{X1}\\
f^{\circ}\left(  v^{\circ}\right)   &  =f^{\bullet}\left(  e^{\mp
\beta\left\vert x_{2}\right\vert \mp\gamma\left\vert x_{1}\right\vert
}v^{\circ}\right)  ,\nonumber
\end{align}
and obtain from the integral identity (\ref{95}) in $\Omega^{\circ}$ the new
one posed in the Sobolev space $H_{0}^{1}\left(  \Omega^{\circ}\right)  $
\begin{align}
a^{\circ}\left(  u^{\circ},v^{\circ}\right)   &  =\left(  \triangledown
u^{\circ},\triangledown v^{\circ}\right)  _{\Omega^{\circ}}\mp\left(  \theta
u^{\circ},\triangledown v^{\circ}\right)  _{\Omega^{\circ}}\pm\left(
\triangledown u^{\circ},\theta v^{\circ}\right)  _{\Omega^{\circ}}-\left(
\theta u^{\circ},\theta v^{\circ}\right)  _{\Omega^{\circ}}-\lambda\left(
u^{\circ},v^{\circ}\right)  _{\Omega^{\circ}}\label{97}\\
&  =f^{\circ}\left(  v^{\circ}\right)  \text{ \ }\forall v^{\circ}\in
H_{0}^{1}\left(  \Omega^{\circ}\right) \nonumber
\end{align}
where $\theta\left(  x\right)  =\left(  \gamma\mathrm{sign~}x_{1}%
,\beta\mathrm{sign~}x_{2}\right)  $. Since $\lambda\notin\sigma^{\circ}$, the
operator $\mathcal{A}_{0,0}^{1\circ}=A^{\circ}$ is an isomorphism and the
problem (\ref{97}) at $\theta=0$ is uniquely solvable in $H_{0}^{1}\left(
\Omega^{\circ}\right)  .$ Furthermore, in view of formulas (\ref{21}) and
(\ref{Cbg}) we have%
\[
\operatorname{Re}a^{\circ}\left(  u^{\circ},v^{\circ}\right)  >\left(
\Lambda^{\star}-\left(  \beta^{2}+\gamma^{2}\right)  -\lambda\right)
\left\Vert u^{\circ};L^{2}\left(  \Omega^{\circ}\right)  \right\Vert ^{2}%
\]
so that the Lax-Milgram Lemma ensures the unique solvability of problem
(\ref{97}), too. The inverse changes (\ref{X1}) give us a solution $u^{1}\in
W_{\pm\beta,\pm\gamma,0}^{1+}\left(  \Omega^{\circ}\right)  $ which falls into
$W_{\pm\beta,\pm\gamma}^{2+}\left(  \Omega^{\circ}\right)  $ due to an
argument in the proof of Lemma \ref{lemmaM2} with a slight modification. We
now multiply $u^{1}$ with the cut-off function $X$ and observe that the
difference $u^{\bullet}-\mathcal{X}u^{1}$ must be find from the problem
(\ref{88}) with the new right-hand side
\begin{equation}
f^{1}=\left(  1-\mathcal{X}^{2}\right)  f^{\bullet}+\left[  \Delta
,\mathcal{X}\right]  u^{1}\in\mathcal{L}_{\pm\beta,\pm\gamma}^{2}\left(
\Omega^{\bullet}\right)  \label{98}%
\end{equation}
which has a support in the strip $\overline{\Xi}^{\circ}\supset\overline{\Xi
}^{+}$. The latter allows us to fix a sufficiently small $\delta>0$ such that
$\mathcal{X}f^{2}\in\mathcal{L}_{\pm\beta+\delta,\pm\gamma}^{2}(\Omega
^{\sharp})$ and the scheme \cite{na17} still works with the weight indexes
$\pm\beta+\delta$ and gives a solution $u^{2}\in W_{\pm\beta+\delta,\pm\gamma
}^{2-}(\Omega^{\sharp})\cap W_{\pm\beta+\delta,\pm\gamma,0}^{1-}%
(\Omega^{\sharp}).$ Multiplying this solution with the cut-off function
(\ref{89}), we observe that, first, $\chi u^{2}\in W_{\pm\beta+\delta
,\pm\gamma}^{2+}(\Omega^{\bullet})\subset W_{\pm\beta,\pm\gamma}^{2+}%
(\Omega^{\circ})$ and, second, it remains to determine the difference
\[
u^{3}=u^{\bullet}-\mathcal{X}u^{1}-\chi u^{2}\in W_{\pm\beta,\pm\gamma}%
^{2+}(\Omega^{\bullet})
\]
from problem (\ref{88}) with the right-hand side%
\begin{equation}
f^{2}=(1-\chi^{2})f^{1}+[\triangle,\chi]u^{2}. \label{99}%
\end{equation}
The last commutator belongs to $W_{\pm\beta+\delta,\pm\gamma}^{1+}%
(\Omega^{\bullet})$ and has a support in the strip $\overline{\Pi^{\sharp
}\left(  1\right)  }$ while the embedding $W_{\pm\beta+\delta,\pm\gamma}%
^{1+}(\Pi^{\sharp}\left(  1\right)  )\subset\mathcal{L}_{\pm\beta,\pm\gamma
}^{2}(\Pi^{\sharp}\left(  1\right)  )$ is compact due to negative increments
of the smoothness and weight exponents. The first term on the right-hand side
of (\ref{99}) has a compact support and a classical construction of a
parametrix in a finite smooth domain gives a compactly supported function
$u^{3}\in H_{0}^{1}(\Omega^{\bullet})$ such that the operator $f^{\bullet
}\rightarrow\mathcal{R}_{\pm\beta,\pm\gamma}^{2}f^{\bullet}=\mathcal{X}%
u^{1}+\chi u^{2}+u^{3}$ gains the necessary properties. We repeat that a
detailed explanation of the above procedure is given in \cite{CaNaTa}. The
operator $\mathcal{A}_{-\beta,-\gamma}^{1\bullet}$ in (\ref{96}) is adjoint
for $\mathcal{A}_{\beta,\gamma}^{1\bullet}$ because the form on the left-hand
side of (\ref{95}) is symmetric. Hence,
\begin{equation}
\mathrm{Ind}\mathcal{A}_{-\beta,-\gamma}^{1\bullet}=-\mathrm{Ind}%
\mathcal{A}_{\beta,\gamma}^{1\bullet}\Longrightarrow\mathrm{Ind}%
\mathcal{A}_{-\beta,-\gamma}^{2}=-\mathrm{Ind}\mathcal{A}_{\beta,\gamma}^{2}.
\label{K1}%
\end{equation}
The implication is supported by the lifting smoothness procedure in Lemma
\ref{lemmaM2}. Moreover, the theorem on the index increment, cf. \cite[\S 3.3
and \S 5.1]{NaPl} ensures the equality%
\begin{equation}
\mathrm{Ind}\mathcal{A}_{-\beta,-\gamma}^{2}=-\mathrm{Ind}\mathcal{A}%
_{\beta,\gamma}^{2}+2 \label{K2}%
\end{equation}
where $2$ is nothing but the total multiplicity of the spectrum of the pencil
$\mathfrak{A}^{\sharp}(\cdot;\lambda)$ in the rectangle $\{\zeta\in
\mathbb{C}:\operatorname{Re}\zeta\in(-\pi,\pi],\ \left\vert \operatorname{Im}%
\zeta\right\vert <\beta\}$, see formula (\ref{39}) and recall our choice of
the upper bound $M^{\sharp}$ in Section \ref{sect2}. Combining (\ref{K1}) and
(\ref{K2}) leads to (\ref{94}). $\boxtimes$

\bigskip

Let us prove the main result of our paper which, owing to the obtained
results, can be obtained in a standard way, see, e.g., \cite[Ch.5]{NaPl}.

\begin{theorem}
\label{TheoremRAD} Let (\ref{Cbg}) and (\ref{888}) be met. Problem (\ref{88})
with the Mandelstam radiation condition has a solution%
\begin{equation}
u\left(  x\right)  =\chi\left(  x_{1}\right)  a^{+}w^{+}\left(  x\right)
+\widetilde{u}\left(  x\right)  \label{K3}%
\end{equation}
with $a^{+}\in\mathbb{C}$ and $\widetilde{u}\in W_{\beta,\gamma}^{2+}%
(\Omega^{\bullet})\cap W_{\beta,\gamma,0}^{1+}(\Omega^{\bullet})$ if and only
if the right-hand side $f$ satisfies the compatibility conditions%
\begin{equation}
(f,v)_{\Omega^{\bullet}}=0\text{ \ }\forall v\in\ker\mathcal{A}_{\beta,\gamma
}^{2}. \label{K4}%
\end{equation}
This solution is defined up to a trapped mode in $\ker\mathcal{A}%
_{\beta,\gamma}^{2}$ with the exponential decay in all directions. The
solution satisfying the orthogonality conditions%
\begin{equation}
(u,v)_{\Omega^{\bullet}}=0\text{ \ }\forall v\in\ker\mathcal{A}_{\beta,\gamma
}^{2} \label{K5}%
\end{equation}
becomes unique and enjoys the estimate%
\begin{equation}
\left\vert a^{+}\right\vert +||\widetilde{u};W_{\beta,\gamma}^{2+}%
(\Omega^{\bullet})||\leq c_{\beta\gamma}(\lambda)\left\Vert f;L_{\beta,\gamma
}^{2}(\Omega^{\bullet})\right\Vert . \label{K6}%
\end{equation}

\end{theorem}

\textbf{Proof.} Owing to Theorems \ref{TheoremInd} and \ref{TheoremASY},
formulas (\ref{K4}) and (\ref{888}) provide a solution $u\in W_{-\beta
,-\gamma}^{2+}(\Omega^{\bullet})\cap W_{-\beta,-\gamma,0}^{1+}(\Omega
^{\bullet})$ together with the representation (\ref{91}) where we need to
eliminate the coefficient $a^{-},$ cf. Section \ref{sect5.2}. To this end, we
observe that $\dim(\ker\mathcal{A}_{-\beta,-\gamma}^{2}\ominus\ker
\mathcal{A}_{\beta,\gamma}^{2})=1$ and there exists a solution $z$ of the
homogeneous problem (\ref{5}), (\ref{6}) in the form%
\begin{equation}
z\left(  x\right)  =\chi\left(  x_{1}\right)  (w^{-}\left(  x\right)
+sw^{+}\left(  x\right)  )+\widetilde{z}\left(  x\right)  \label{K7}%
\end{equation}
where $\widetilde{z}\in W_{\beta,\gamma}^{2+}(\Omega^{\bullet})$ and
$s\in\mathbb{C}$ is the reflexion coefficient, $\left\vert s\right\vert =1.$
Finally, the difference $u-a^{-}\zeta$ satisfies the radiation condition and
takes the form (\ref{K3}). A solution $z$ in $\ker\mathcal{A}_{-\beta,-\gamma
}^{2}\setminus\ker\mathcal{A}_{\beta,\gamma}^{2}\neq\varnothing$ has at least
one non-trivial coefficient $a^{\pm}$ in its representation (\ref{91}). Let us
assume that $a^{+}=0$ and, therefore, representation (\ref{K7}) is not
possible. We truncate the domain $\Omega^{\bullet}$ like $\Omega_{R}^{\bullet
}=\{x\in\Omega^{\bullet}:\left\vert x_{j}\right\vert <R,$ $\ j=1,2\}.$ We also
denote $T_{R}^{\bullet}=\{x\in\Omega^{\bullet}:x_{1}=R\}$ and insert $z$ into
the Green formula on $\Omega_{R}^{\bullet}.$ Taking the Dirichlet condition
(\ref{6}) into account, we have%
\begin{equation}
0=\int_{\partial\Omega_{R}^{\bullet}\setminus\partial\Omega^{\bullet}%
}(\overline{z\left(  x\right)  }\partial_{n}z\left(  x\right)  -z\left(
x\right)  \overline{\partial_{n}z\left(  x\right)  })ds_{x} \label{K8}%
\end{equation}
where $\partial_{n}$ is the outward normal derivative. The exponential decay
of $\widetilde{z}\left(  x\right)  $ as $\left\vert x\right\vert
\rightarrow+\infty$ allows us to get rid in (\ref{K8}) of the remainder
$\widetilde{z}$ and the whole integral over $(\partial\Omega_{R}^{\bullet
}\setminus\partial\Omega^{\bullet})\setminus T_{R}^{\bullet}.$ Moreover, we
add the integrals along $\{x\in\Omega^{\bullet}:x_{1}=R,$ $\ \pm x_{2}>R\}$
and write%
\[
\left\vert a^{+}\right\vert ^{2}\int_{T_{R}^{\bullet}}\left(  \overline
{w^{+}(x)}\frac{\partial w^{+}}{\partial x_{1}}(x)-w^{+}(x)\overline
{\frac{\partial w^{+}}{\partial x_{1}}(x)}\right)  dx_{2}=O(e^{-\min
\{\beta,\gamma\}R}).
\]
Finally, we integrate in $R\in(N,N+1)$ and send $N\in\mathbb{N}$ to infinity
to obtain $\left\vert a^{+}\right\vert ^{2}q(w^{+},w^{+})=0$ according to
(\ref{F2}). Recalling (\ref{F3}) and (\ref{F5}), (\ref{FFF}), we conclude that
$a^{+}=0$, $z\in\ker\mathcal{A}_{\beta,\gamma}^{2}$ and come across a
contradiction to our assumption. The equality $\left\vert s\right\vert =1$ for
the reflexion coefficient is verified by a similar calculation based on
bi-orthogonality conditions of type (\ref{F5}). $\boxtimes$

\subsection{Available generalizations\label{sect5.4}}

Many unnecesary restrictions were introduced in our paper to simplify
demonstration only. In particular, the Laplace operator $\Delta$ can be
replaced by a formally self-adjoint second-order differential operator in the
divergence form $\nabla\cdot\mathbf{A}(x)\nabla$ with a positive definite
symmetric $2\times2$-matrix $\mathbf{A}$ with periodic measurable bounded
coefficients. The boundary $\partial\varpi$ of the periodicity cell can be,
e.g., Lipschitz while the semi-infinite foreign inclusion can be formed by
varying the boundary and the coefficents.

\begin{figure}
\begin{center}
\includegraphics[scale=0.45]{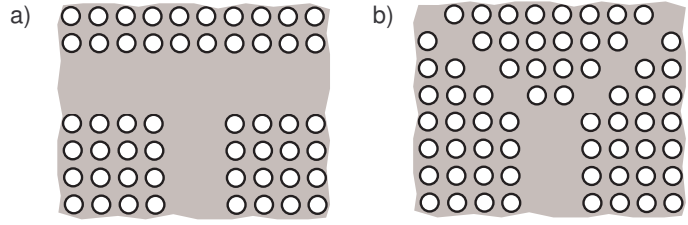}
\end{center}
\caption{Open waveguides in the double-periodic plane with several outlets to
infinity.}
\label{fig.5}
\end{figure}

As in \cite{CaNaTa}, the open waveguide may have several outlets to infinity,
cf. fig. \ref{fig.5}.%

Each of open waveguides may enjoy a local perturbation, cf. fig. \ref{fig.2},
a and b. By means of the classical approach \cite{Jones} one can readily
detect a point of the discrete spectrum. Let us show the existence of at least
one eigenvalue $\lambda^{\blacksquare}\in(0,M_{1}(0))$ in the spectrum of the
operator $A^{\blacksquare}$ of the Dirichlet problem (\ref{5}), (\ref{6}) in
the domain $\Omega^{\blacksquare}=\Omega^{\circ}\cup\Xi^{+}\cup\Xi
^{\blacksquare}$ where $\Xi^{\blacksquare}$ is the rectangle $(-J_{1}%
l_{1},J_{1}l_{1})\times(-J_{2}l_{2},J_{2}l_{2})$ where $J_{1},$ $J_{2}%
\in\mathbb{N}$, see fig. \ref{fig.2}, a. The lower bound of the spectrum
$\sigma^{\blacksquare}$ of $A^{\blacksquare}$ is still equal to $M_{1}(0).$ We
choose $J_{1}$ and $J_{2}$ such that the Dirichlet problem in $\Xi
^{\blacksquare}$ has the principal eigenpair
\[
\lambda^{\blacksquare}=\frac{\pi^{2}}{4}\left(  \frac{1}{J_{1}^{2}}+\frac
{1}{J_{2}^{2}}\right)  <M_{1}(0),\text{ \ }u^{\blacksquare}(x)=\cos\left(
\frac{\pi}{J_{1}}x_{1}\right)  \cos\left(  \frac{\pi}{J_{2}}x_{2}\right)
\]
Extending $u^{\blacksquare}$ as null from $\Xi^{\blacksquare}$ onto
$\Omega^{\blacksquare}$, we apply the min principle, cf. \cite[Thm
10.2.1]{BiSo}, and obtain
\[
\underline{\sigma^{\blacksquare}}=\underset{u\in H_{0}^{1}(\Omega
^{\blacksquare})\setminus\{0\}}{\inf}\frac{\left\Vert \nabla u;L^{2}%
(\Omega^{\blacksquare})\right\Vert ^{2}}{\left\Vert u;L^{2}(\Omega
^{\blacksquare})\right\Vert ^{2}}\leq\frac{\left\Vert \nabla u^{\blacksquare
};L^{2}(\Xi^{\blacksquare})\right\Vert ^{2}}{\left\Vert u^{\blacksquare}%
;L^{2}(\Xi^{\blacksquare})\right\Vert ^{2}}=\lambda^{\blacksquare}<M_{1}(0).
\]
Hence, the point $\underline{\sigma^{\blacksquare}}$ of the spectrum of the
operator $A^{\blacksquare}$ belongs to its discrete spectrum.

\end{document}